\documentclass[a4paper, 11pt]{article}
\pdfoutput=1
\usepackage{amsmath}
\usepackage{amssymb}
\usepackage{graphicx}

\usepackage{epsfig}
\usepackage[english]{babel}
\usepackage{lineno}
\usepackage{color}

\usepackage{epstopdf}

\newcommand{\be}{\begin{equation}}
\newcommand{\ee}{\end{equation}}
\newcommand{\bea}{\begin{eqnarray}}
\newcommand{\eea}{\end{eqnarray}}

\newcommand{\nod}{\noindent}

\newcommand{\ba}{\begin{array}}
\newcommand{\ea}{\end{array}}
\newcommand{\bc}{\begin{center}}
\newcommand{\ec}{\end{center}}

\DeclareMathOperator{\Tr}{Tr}
\DeclareMathOperator{\sign}{sign}

\newtheorem{prop}{Proposition}

\begin{document}

\title{Oscillating epidemics in a dynamic network model: stochastic and mean-field analysis}

\author{Andr\'as Szab\'o-Solticzky $^1$, Luc Berthouze $^{2}$, Istvan Z. Kiss $^{3,\ast}$ \& P\'eter L. Simon $^{1}$}

\maketitle

\begin{center}
$^{1}$ Institute of Mathematics, E\"otv\"os Lor\'and University Budapest, and \\ Numerical Analysis and Large Networks Research Group, Hungarian Academy of Sciences, Hungary\\
$^{2}$ Centre for Computational Neuroscience and Robotics, University of Sussex, Falmer, Brighton BN1 9QH, UK\\
$^{3}$ School of Mathematical and Physical Sciences, Department of Mathematics, University of Sussex, Falmer, Brighton BN1 9QH, UK
\end{center}

\vspace{1cm}

\begin{abstract}
An adaptive network model using $SIS$ epidemic propagation with link-type dependent link activation and deletion is considered. Bifurcation analysis of the pairwise ODE approximation and the network-based stochastic simulation is carried out, showing that three typical behaviours may occur; namely, oscillations can be observed besides disease-free or endemic steady states. The oscillatory behaviour in the
stochastic simulations is studied using Fourier analysis, as well as through analysing the exact master equations of the stochastic model. A compact pairwise approximation for the dynamic network
case is also developed and, for the case of link-type independent rewiring, the outcome of epidemics and changes in network structure are concurrently presented in a single bifurcation diagram.
By going beyond simply comparing simulation results to mean-field models, our approach yields deeper insights into the observed phenomena and help better understand and map out the limitations of mean-field models.
\end{abstract}

\nod {\bf Keywords:} SIS epidemic; pairwise model; dynamic network; oscillation \\

\vspace{1cm}


\vspace{1cm}
\begin{flushleft}
$\ast$ corresponding author\\
email: i.z.kiss@sussex.ac.uk\\
\end{flushleft}

\newpage

\section{Introduction}

Network-related research has recently seen rapid developments in the domain of dynamic or adaptive networks \cite{ThiloAdCoevNetw, SaldanaJMB, MarceauPRE,
ShawSchwartz}. This is partly motivated by strong empirical observations that in many instances dynamical processes on networks co-evolve with the dynamics of the networks themselves
\cite{ThiloAdCoevNetw}. This in turn leads to a wider spectrum of possible behaviours, such as bistability and oscillations, when compared to static network models. Examples of dynamic/adaptive
networks
are abundant and it is a common feature of both technological and social networks. In this paper, we study dynamic networks in which the timescales of the two processes are comparable, and we do
this in the context of a basic but fundamental model of disease transmission.

Early work in the area of dynamic networks, which then gave rise to many model improvements and extensions, concentrated on epidemic models with `smart' rewiring, where infection transmitting
links are replaced by non-risky ones, with a link-conserving network dynamics. The most widely used approach to study such systems is the development of mean-field models, such as pairwise
\cite{ThiloPhyRevLett,
PRAcikk} and effective-degree models \cite{MarceauPRE, TaylorPhysRevE}, which usually manage to capture and characterise such processes to a good level of detail. For example, using pairwise
models and simulations, Gross et al. \cite{ThiloPhyRevLett} showed that rapid rewiring can stop
disease transmission, while Marceau et al. \cite{MarceauPRE} used the effective-degree model, i.e. an improved compartmental model formalism, to obtain better predictions and shed more light on the
evolution of the network structure. Furthermore, Juher et al. \cite{SaldanaJMB} fine-tuned the bifurcation analysis of the model originally proposed in \cite{ThiloPhyRevLett} to distinguish between two
extinction scenarios and to provide an analytical
condition for the occurrence of a bistability region.

Recognising the idealised nature of link-number-conserving ÂsmartÂ rewiring, various relaxations to these models were introduced. For example, Kiss et al. \cite{PRAcikk} and Szabo et al.
\cite{Szabo2012DEA} proposed a number of dynamic network models including random link activation-deletion with and without constraints, as well a model that considered non-link preserving
link-type-dependent activation deletion. Taylor et al. \cite{TaylorPhysRevE} used the effective degree formalism to analyse a similar random link activation deletion model and to highlight the potential power of link
deletion in eradicating epidemics. Recently Rogers et al. \cite{Rogers2012} introduced an SIRS
model where random link activation deletion is combined with `smart' rewiring. They studied the resulting stochastic model at the level of singles and pairs as well as its deterministic limit, i.e. pairwise or pair-based models. Other node dynamics, such as the voter model, have also been extensively studied on adaptive networks; see \cite{ThiloPhysD2014, ThiloCoMaAppl2013} for a review.

An overview of the above studies highlights some important modelling and analysis challenges. On the model development front it is clear that many models are still very idealised and small steps toward
increasing model realism can quickly lead to a disproportionate increase in model complexity. The agreement between mean-field and stochastic models is very much model
and parameter dependent with potentially large parameter regions in which agreement can be either good or poor. This is especially the case when oscillations are encountered or
expected \cite{GrossKevrekidis}. When oscillatory behaviour is expected, there are only few results in an epidemic context, see \cite{ZhouetalPRE2012} for example. Studies of epidemic
propagation on adaptive networks have not focused much on characterising oscillatory behaviour in simulations and there is no widely accepted analytical framework for this as of yet.  Another major drawback in
the analysis of dynamic networks is that even when mean-field models perform well, these only provide limited insight about the structure and evolution of the network.

In this paper, we set out to address some of these issues and make important steps towards a more satisfactory treatment and analysis of dynamic network models, along with a better use and
integration of various mathematical methods and techniques that can be employed. To this end we propose a dynamic network model using $SIS$ epidemic propagation with both random and
preferential, e.g., link-type based, link activation and deletion. We provide model
analysis based on both (a) the mean-field and (b) the purely network-based stochastic simulation model. The latter is carried out using Fourier analysis, as well as through analysing the exact stochastic
model in terms of the master equations. Further, we develop a compact pairwise approximation for the current dynamic network model and, for the case of random activation deletion, we provide a
bifurcation analysis of the network structure itself, aiming for a more comprehensive analysis, whereby both system-level (e.g., the outcome of the epidemic) and network-level (e.g., the achievable network characteristics of
network types) behaviours are concurrently analysed and characterised.

The paper is structured as follows. The link-type-dependent activation-deletion model is formulated in Section \ref{ModelForm}. Then the pairwise ODE approximation augmented by terms accounting for
preferential link cutting and creation is analysed in Section \ref{sec:PW}.  Guided by the bifurcation analysis of the mean-field model, in Section \ref{sec:NetSim}, we carry out a detailed study of the
agreement between models and simulations. Oscillations are observed both in the mean-field approximation and the stochastic network simulation. In Section \ref{sec:ExplOsc}, the emergence of oscillations in
the stochastic model is investigated based on the outcome of simulation and on the exact master equations. Finally, in Section \ref{sec:NetBif}, a detailed and common bifurcation map of the system and
network behaviour is presented emphasising its benefits in achieving a fuller understanding of the model as a whole, especially for dynamic network models.


\section{Model formulation} \label{ModelForm}

In this paper $SIS$ (susceptible-infectious/infected-susceptible) epidemic propagation is considered on an adaptive network with link-type dependent link activation and deletion. Specifically, the model
incorporates the following independent Poisson processes:
\begin{itemize}
\item {\bf Infection:} Infection is transmitted across each contact between an $S$ and an $I$ node, or $(SI)$ link, at rate $\tau$,
\item {\bf Recovery:}  Each $I$ node recovers at rate $\gamma$, and this is independent of the network,
\item {\bf Link activation:} A non-existing link between a node of type $A$ and another of type $B$ is activated at rate $\alpha_{AB}$, with $A,B \in \{S,I\}$,
\item {\bf Link deletion:} An existing link between a node of type $A$ and another of type $B$ is terminated at rate $\omega_{AB}$, with $A,B \in \{S,I\}$.
\end{itemize}
We note that once a link type is chosen, the activation or deletion of such a link is done at random. Our model is significantly different from the widely used setup of `smart' rewiring
\cite{ThiloAdCoevNetw}, where the $S$ nodes have full knowledge of the states of all other nodes and choose to minimise their exposure to the epidemic by cutting links to $I$ neighbours and
immediately rewiring to a randomly chosen $S$ node. This also conserves the number of links in the network and can make the analysis more tractable, by reducing the complexity of the model.

Here, we set out to explore and explain as fully as possible the complete spectrum of system behaviours, including classical bifurcation analysis at system level, e.g., die-out, endemic equilibria or
oscillations, as well as the evolution of the network structure and attainable network equilibria. In order to do this, we will employ a number of approaches including: (a) an exact Master Equation
formulation for small networks, (b) full network-based Monte Carlo simulations and (c) two different types of pairwise or pair-based mean-field ODE models. As regards the rewiring parameters we focus on two scenarios, namely:
\begin{enumerate}
\item[\textbf{A.}] $\alpha_{SI}=\alpha_{II}=0$ and $\alpha_{SS} \neq 0$, and $\omega_{II}=\omega_{SS}=0$ and $\omega_{SI} \neq 0$, and
\item[\textbf{B.}] $\alpha_{SI}=\alpha_{II}=\alpha_{SS}=\alpha$ and $\omega_{SI}=\omega_{II}=\omega_{SS}=\omega$.
\end{enumerate}
While the first is motivated by practical considerations, such as those used in the `smart' rewiring -- where nodes aim to minimise the risk of becoming infected while maintaining their connectivity to the
network, the second scenario removes the dependency of link activation and deletion on pair type and leads to a simple, more tractable model.

\section{Simple pairwise model: bifurcation analysis} \label{sec:PW}

We start by formulating the pairwise model for the expected values of the node and pair numbers.
As was shown in \cite{PRAcikk}, this gives rise to
\begin{align}
\dot{[I]}&=\tau [SI] - \gamma [I],  \label{SMF1} \\
\dot{[SI]}&=\gamma ( [II] - [SI])+ \tau ([SSI]-[ISI]-[SI])+\alpha_{SI} ([S][I]-[SI])-\omega_{SI}[SI],  \label{SMF2} \\
\dot{[II]}&=-2\gamma [II] +2\tau ([ISI]+[SI])+\alpha_{II} ([I]([I]-1)-[II])-\omega_{II}[II],  \label{SMF3} \\
\dot{[SS]}&=2\gamma [SI] - 2\tau [SSI]+\alpha_{SS} ([S]([S]-1)-[SS])-\omega_{SS}[SS].  \label{SMF4}
\end{align}
From the model it follows that $[S]+[I]=N$, and that the closures of the triples in terms of singles and pairs require the expected average degree of $S$ nodes.
This is given by
\begin{equation}
k_{S}(t)=\frac{[SS]+[SI]}{[S]}. \label{ks}
\end{equation}
The well-known closures \cite{Keeling1999} are used, namely
\begin{equation}
[SSI]=\frac{(k_{S}-1)[SS][SI]}{k_S[S]} \,\,\, \text{and} \,\,\, [ISI]=\frac{(k_{S}-1)[SI][SI]}{k_S[S]}.
\end{equation}
Upon applying these closures, a self-consistent system with 4 ODEs is obtained. This can be analysed using
classical bifurcation theory techniques. Focusing on scenario \textbf{A}, the system admits two equilibria: (a) a disease-free equilibrium $ ([S], [I], [SI], [II], [SS])=(N, 0, 0, 0,
N(N-1))$ and (b) an endemic equilibrium which cannot be explicitly given due to being the solution of a quartic equation.

\begin{figure}[h!]
  \centering
	\includegraphics[width=0.6\textwidth, height=0.3\textheight]{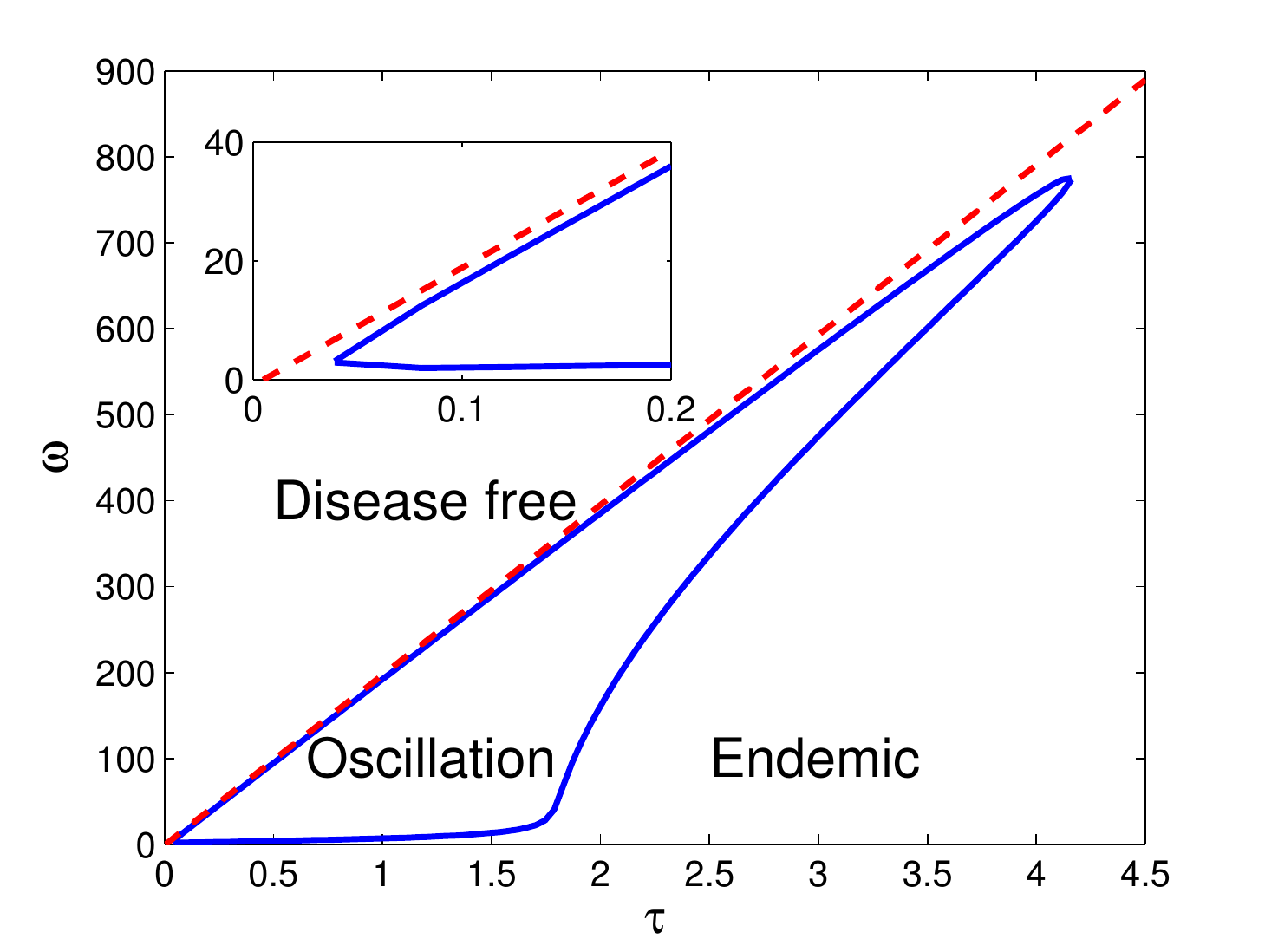}
	\caption{Bifurcation diagram for the pairwise ODE model in the $(\tau,\omega_{SI})$ parameter space for $N=200$, $\gamma=1$ and $\alpha_{SS}=0.04$. The transcritical bifurcation occurs along the dashed line, and the Hopf bifurcation occurs along the perimeter of the island.} \label{SMF_Bif}
\end{figure}

The linerisation around the disease-free steady state gives rise to a $4 \times 4$ Jacobian, the eigenvalues of which can be determined explicitly, see Appendix A.
As shown, two of the eigenvalues are always negative and the remaining two have negative real part if and only if
\begin{equation}
\omega_{SI}>\tau (N-2) - \gamma,
\end{equation}
which gives rise to a transcritical bifurcation where $\omega_{SI}=\tau (N-2) - \gamma$, see Fig.~\ref{SMF_Bif}. Thus the following Proposition holds.
\begin{prop}
The disease-free steady state is stable if and only if $ \omega_{SI} >\tau (N-2)-\gamma $. \label{prop1}
\end{prop}

\begin{figure}[h!]
	\centering
       		\includegraphics[scale=0.45]{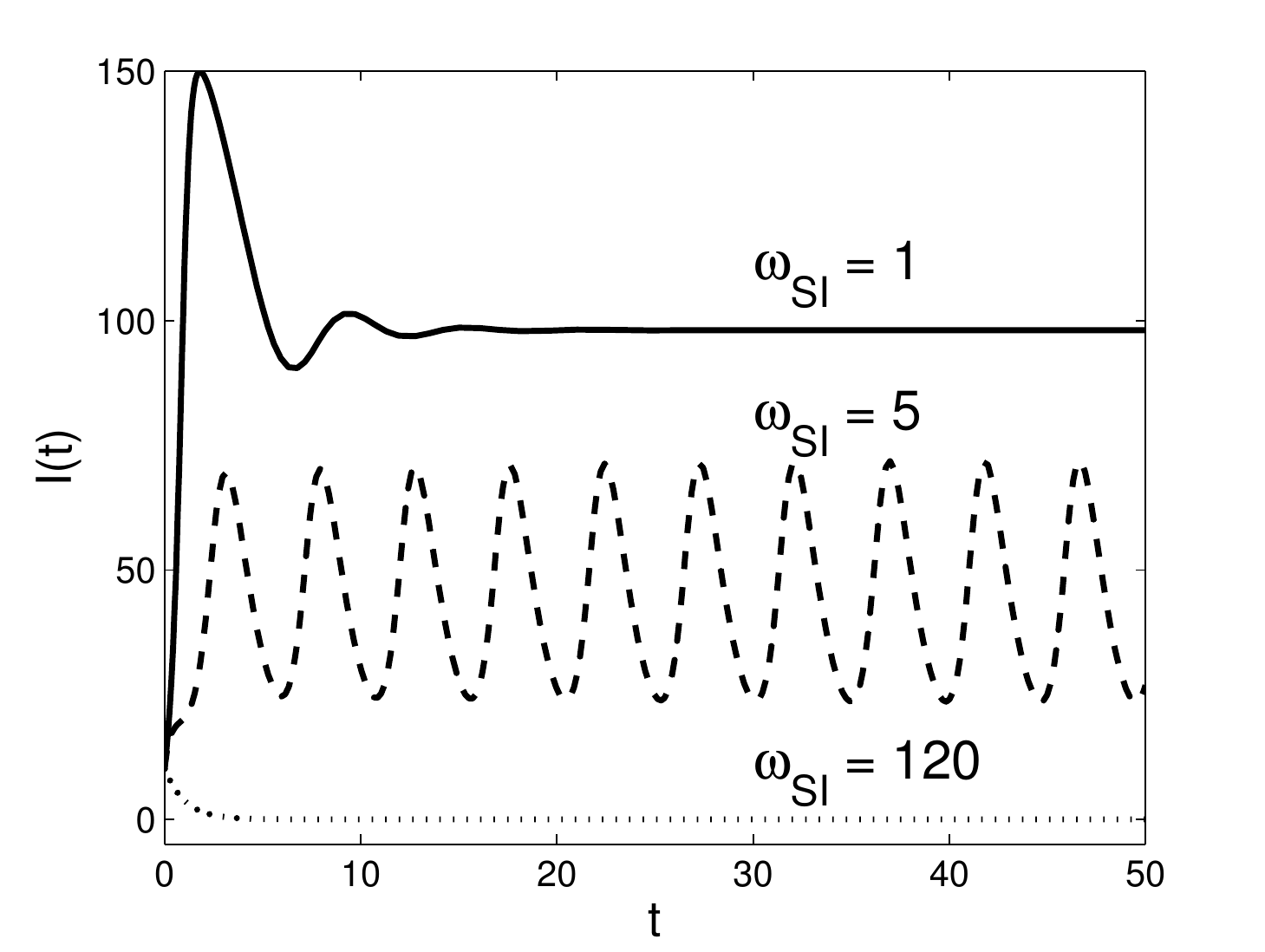}
                	\includegraphics[scale=0.45]{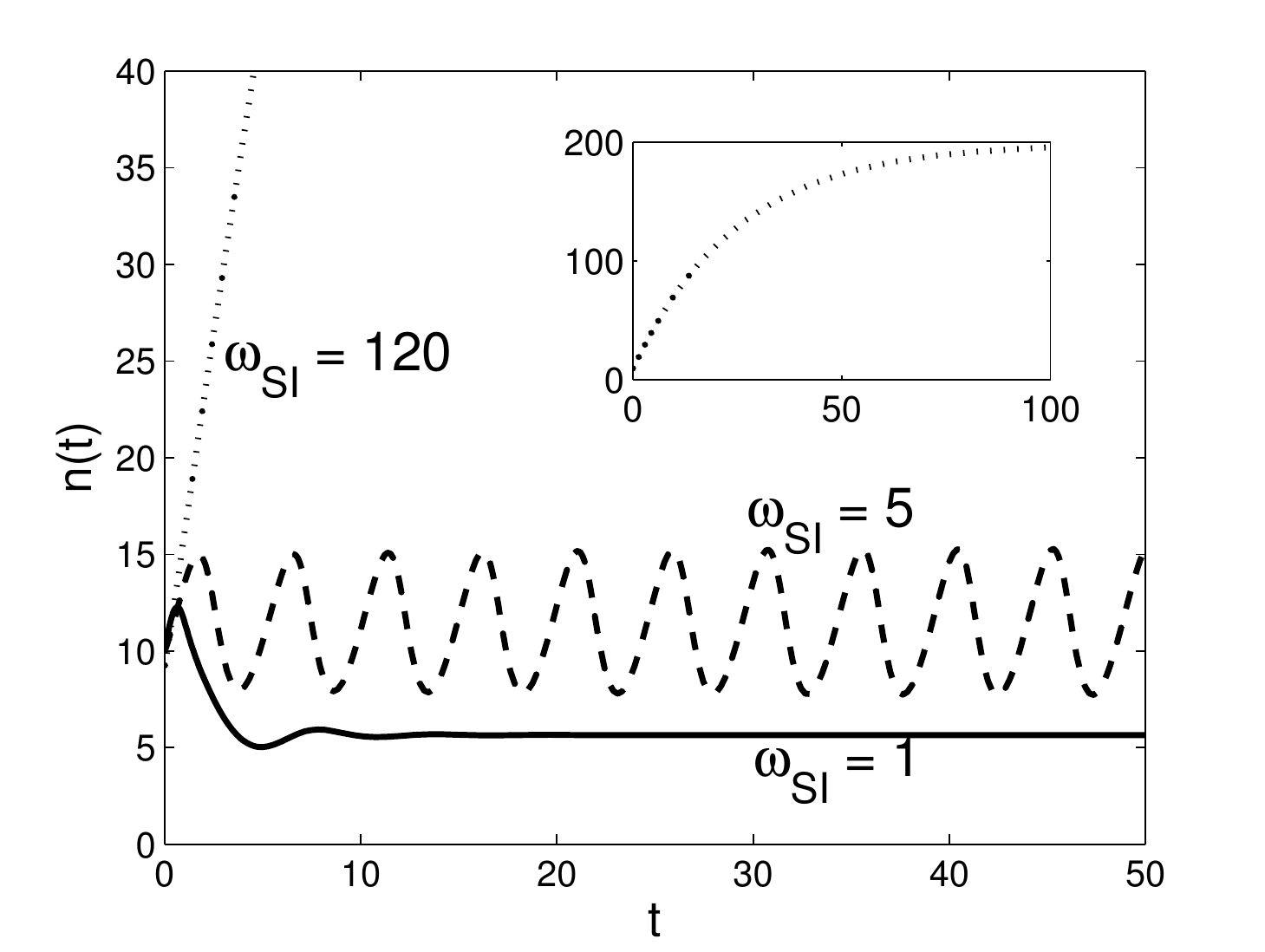}
        \caption{Time dependence of prevalence (left panel) and average degree (right panel) in the pairwise ODE model for $\omega_{SI}=1$ (continuous), $\omega_{SI}=5$ (dashed) and for $\omega_{SI}=120$ (dotted). The inset in the right panel shows the time dependence of the average degree for $\omega_{SI}=120$. The values of the other parameters are fixed at $\tau =0.5$, $N=200$, $\gamma=1$ and $\alpha_{SS}=0.04$. }\label{SMF_TimeEvol}
\end{figure}

As mentioned above the endemic steady state is the solution of a quartic equation, see Appendix A for its detailed derivation.
The analysis of this equation leads to the following proposition concerning the existence of the endemic steady state.
\begin{prop}
If $x\in (0,N)$ is a root of polynomial \eqref{quartic} given in Appendix A, then the system has an endemic steady state, the coordinates of which can be given as
$$
[S]_{ss} = x, \quad [I]_{ss} =N-x, \quad [SI]_{ss} = \frac{\gamma}{\tau}(N-x),
$$
$$
[SS]_{ss} = x(x-1)-2\frac{\omega_{SI}\gamma}{\alpha_{SS}\tau} (N-x), \quad [II]_{ss} = \frac{\gamma (N-x)^2}{\tau x} +\frac{(N-x)[SS]_{ss}}{[SS]_{ss}+[SI]_{ss}} .
$$
 \label{propendemic}
\end{prop}
An extensive numerical study has shown that below the line of transcritical bifurcation there is a unique endemic steady state, i.e., the polynomial has a single root providing a biologically plausible steady
state, where the values of all singles and pairs are positive. For the endemic steady state, the coefficients of the characteristic polynomial can only be determined numerically but this does not prevent from
determining where the Hopf bifurcation arises (see Appendix A for details). The Hopf bifurcation points form a set as depicted by the curve, i.e., the perimeter of the island, in Fig.~\ref{SMF_Bif}. The results
of the numerical study show that within the Hopf island the mean-field model exhibits stable oscillations. The region below the transcritical bifurcation line and outside the Hopf island is where the
endemic
equilibrium is stable. It is important to note that the system-level analysis can be complemented by the observation that the expected average degree displays a behaviour similar to that of the expected
number of infected, as illustrated by Fig.~\ref{SMF_TimeEvol}.

A careful analysis of the plots based on the mean-field model allows to make the following important remarks.
Edge addition and creation acts on potentially very different scales or number of edges.
The creation of $(SS)$ links acts on a set of edges with cardinality of $O(N^2)$, while the infection process and the deletion of $(SI)$ edges act on a set of edges whose number scales as $\langle k
\rangle [I]$. In order to ensure that all processes act in a comparable way, i.e., on the same timescale, rates need to be adjusted accordingly, e.g., smaller values for the $\alpha_{SS}$ rate compared to
the rate of cutting $SI$ links, $\omega_{SI}$. The conclusion of the analysis of the pairwise model is that the system can exhibit three different behaviours: (a) disease-free steady state, (b) endemic
steady state, and (c) oscillations. Moreover, the bifurcation boundaries separating these states can be analytically determined and for given parameter values numerically computed. This analysis offers a
valuable insight into the possible behaviours that can be expected from the full network simulation.

Turning to scenario \textbf{B}, i.e., the case where $\alpha_{SS}=\alpha_{SI}=\alpha_{II}:=\alpha$ and $\omega_{SS}=\omega_{SI}=\omega_{II}:=\omega$, we determine the steady states and the local
behaviour around them. In this case the coordinates of the disease-free steady state are $[I]=0$, $[SI]=0$, $[II]=0$ and $ [SS]=\frac{N(N-1)\alpha}{\omega+\alpha}$. The Jacobian matrix corresponding
to
this steady state is
\begin{equation}
J = \left( \begin{array}{cccc}
-\gamma & \tau & 0 & 0  \\
\alpha N & -\gamma -\tau -\omega - \alpha +Q & \gamma & 0  \\
-\alpha & 2\tau & -2\gamma-\omega-\alpha & 0  \\
\alpha(-2N+1) & 2\gamma -2Q & 0 & -\omega-\alpha
\end{array} \notag
\right),
 \label{jacobi}
\end{equation}
where
$$
Q = \frac{\tau (n-1) (N-1)\alpha}{n(\omega+\alpha)} .
$$
Solving the equation $\det J = 0$ for $\tau$ yields that the transcritical bifurcation occurs at
\begin{equation}
\tau_c = \frac{n\gamma(\omega+\alpha)(2\gamma+\omega+\alpha)(\omega+\alpha+\gamma)}{n(\omega+\alpha)(\alpha N - 2\gamma\alpha - \gamma\omega) + (n-1)(N-1)\alpha\gamma } \ . \label{tauTC}
\end{equation}
Numerical investigation shows that for $\tau <\tau_c$ the solutions of the system tend to the disease-free steady state, while for $\tau >\tau_c$ the solutions converge to the endemic steady state.
Oscillations were not observed in this case.

\section{Characterisation of the full network-based stochastic model} \label{sec:NetSim}

With the insight gained from the analysis of the mean-field model, we set out to investigate and map out the
spectrum of behaviour directly from the network-based stochastic simulation, that is, the analogue of the continuous-time Markov Chain corresponding to the
process. The simulation is based on a careful and continuous book keeping of all possible single events, i.e., infection across a link, recovery of a node,
activation or deletion of a link. Based on a 'current' state all possible events and their rates are computed, or these are known based on a previous state supplemented by necessary changes induced by
the most recent event. The time to the next event is chosen from an exponential distribution whose parameter is the total rate. Then, an event is chosen at random but proportionally to its rate. Given that all existing and non-existing links need to be accounted for, the algorithm which includes the storage, update and referencing back and forth between
rates
and events becomes more complex, concretely, from order $N$ to order $N^2$ complexity.

Simulations were run in which $\tau$ and $\omega$ were systematically sampled in the interval $[0.05,20]$ and $[0.05,13]$ respectively. For each configuration pair, $100$ simulations were run with
$N=200$, $\gamma=1$, $\alpha_{SS}=0.04$, $T_{max}=1320$ and resampling rate $0.01$, yielding 132000 time points per run.

\subsection{Simulation results}

\begin{figure}[h!]
  \centering
	\includegraphics[scale=0.5]{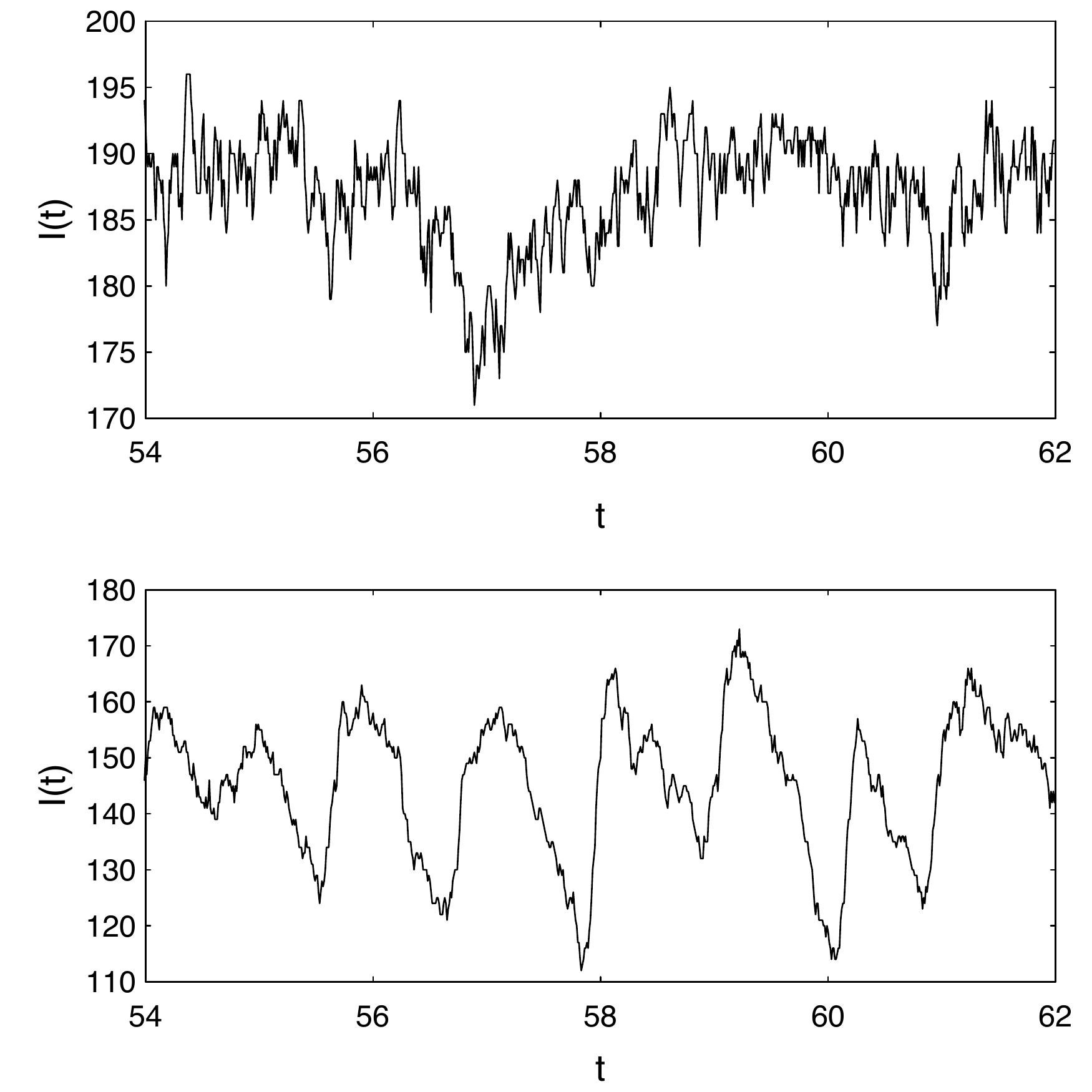}
	\caption{Sample time series ($800$ time points) of the number of infected nodes in the endemic regime (top panel, $\omega_{SI}=0.5$) and oscillatory regime (bottom panel, $\omega_{SI}=4.0$). These samples were randomly chosen from one of $100$ realisations using the following parameters: $N=200$, $\tau=12$, $\gamma=1$, $\alpha_{SS}=0.04$.} \label{fig:trajexamples}
\end{figure}

Simulations show that the system exhibits exactly those three regimes that were predicted by the pair-based mean-field model. Figure~\ref{fig:trajexamples} shows representative samples of the system's
behaviour in the endemic and oscillatory regimes. As we cannot expect sharp bifurcation boundaries between regimes in the stochastic model, we now propose empirical definitions for the boundaries of
the different regimes.

\begin{figure}[h!]
  \centering
	\includegraphics[scale=0.5]{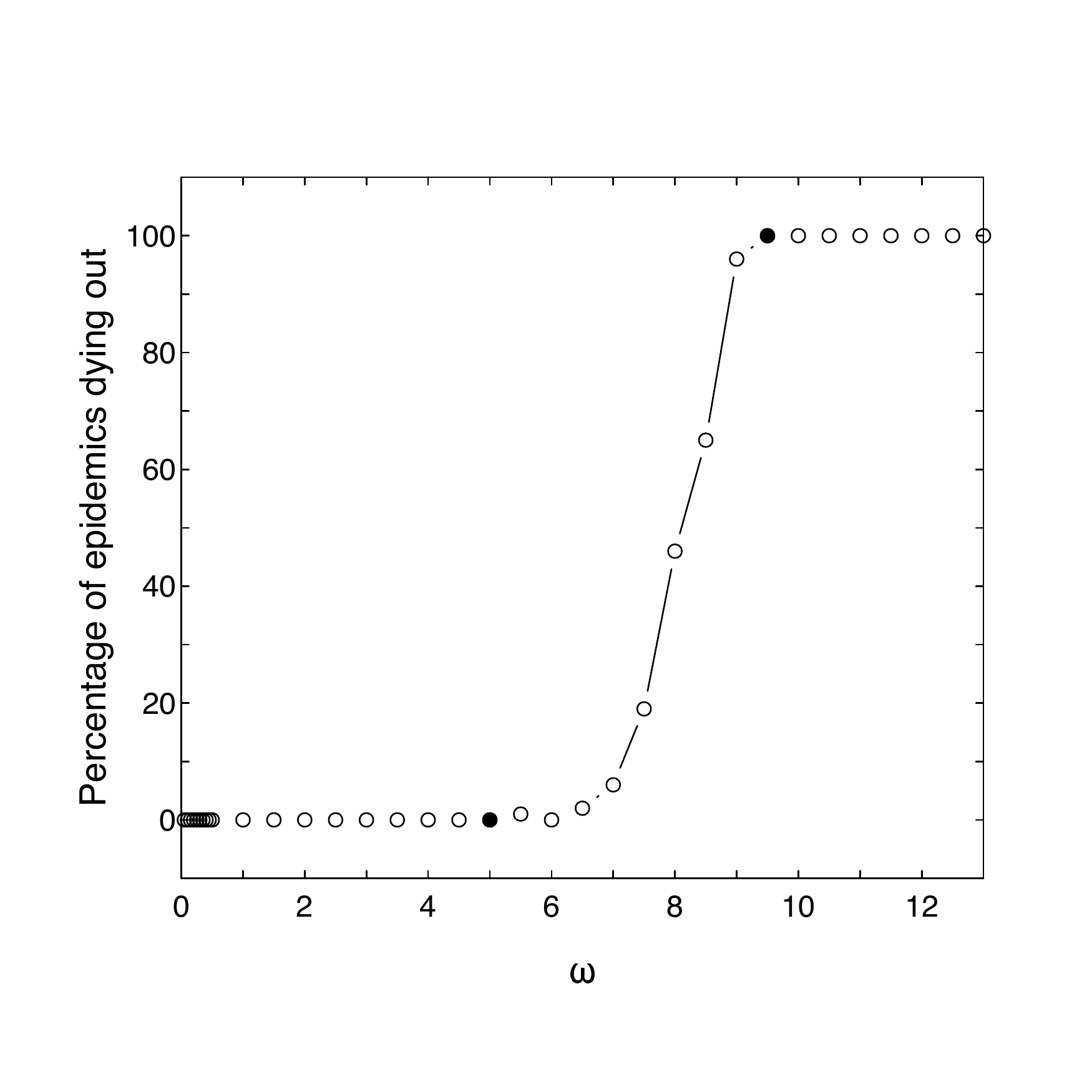}
	\caption{Percentage of epidemics dying out before $T_{max}=1320$ with $\tau=12$ and $\omega_{SI}$ (horizontal axis) varying between $0.05$ and $13$. This percentage is calculated out of $100$ realisations using the following parameters: $N=200$, $\tau=12$, $\gamma=1$, and $\alpha_{SS}=0.04$. The values of $\omega_{SI}$ below which no realisations die out and above which all realisations die out define two possible boundaries for the disease-free regime and are shown by solid circles. } \label{fig:zerodiseaseboundary}
\end{figure}

\textit{Boundary of the disease-free regime.} For each configuration pair $(\tau,\omega)$, we determined the proportion of realizations in which the epidemic died out. Figure~
\ref{fig:zerodiseaseboundary}
shows a representative example (at $\tau=12$) when $\omega$ is varied over the whole range $[0.05,13]$ and all other parameters are kept constant. The figure shows a sharp increase in the proportion of
realisations in which epidemics die out such that two potential boundaries can be defined: one which corresponds to the parameter up to which none of the realisations die out, and one in which all realisations die out. These two boundaries delimit a strip that can be considered as the bifurcation boundary of the disease-free regime. The width of this strip is shown by the grey-shaded area in Figure~\ref{fig:stochbifdiag}.

\begin{figure}[h!]
  \centering
	\includegraphics[scale=0.5]{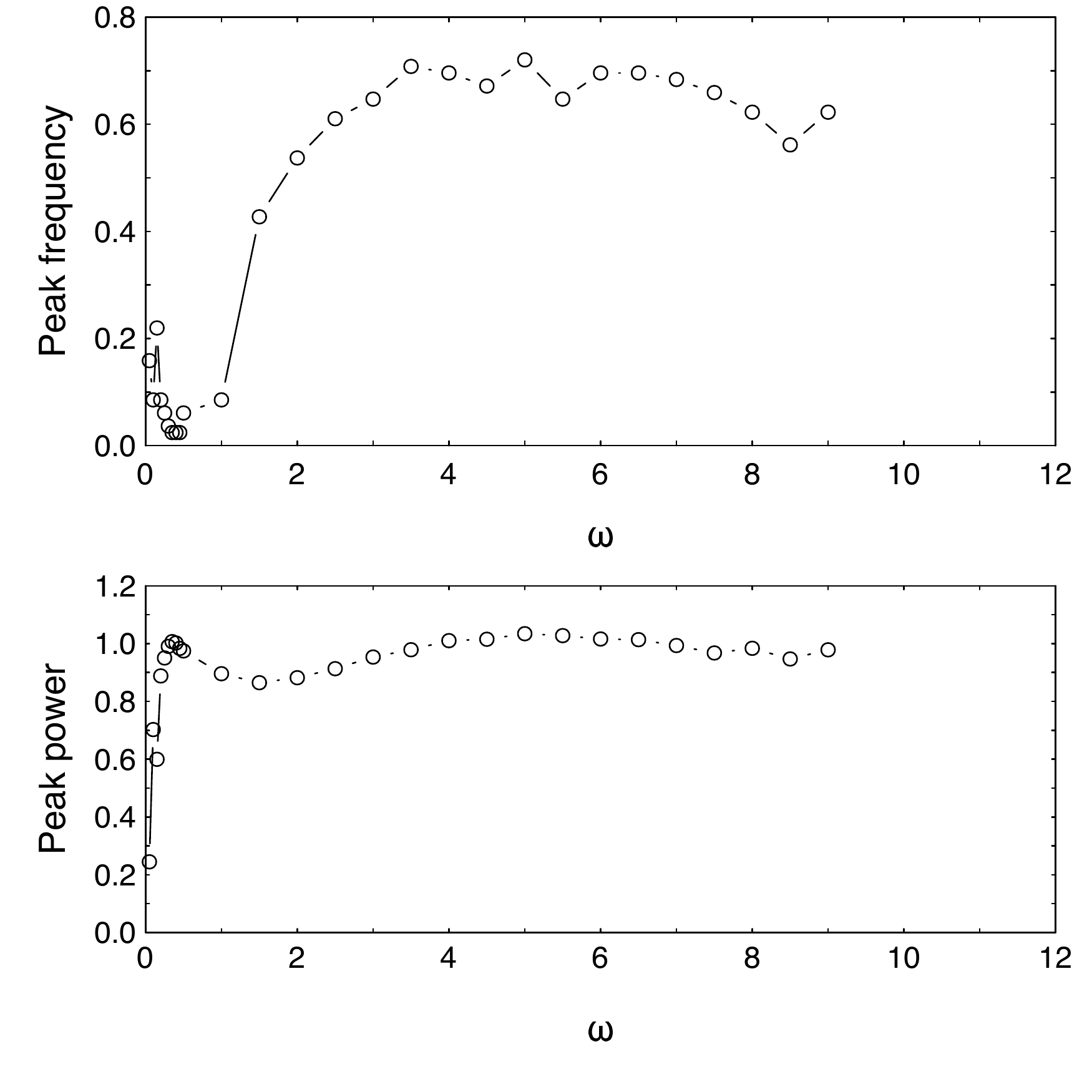}
	\caption{Frequency at peak power and peak power for $\tau=12$ and $\omega_{SI}$ (horizontal axis) varying between $0.05$ and $13$. The absence of data for $\omega_{SI}>9$ reflects the fact that the disease-free regime has been reached, see Fig.~\ref{fig:zerodiseaseboundary}. Thresholding of peak frequency at non-zero value makes it possible to define a boundary for the oscillatory regime. }\label{fig:oscboundary}
\end{figure}

\begin{figure}[h!]
  \centering
	\includegraphics[scale=0.5]{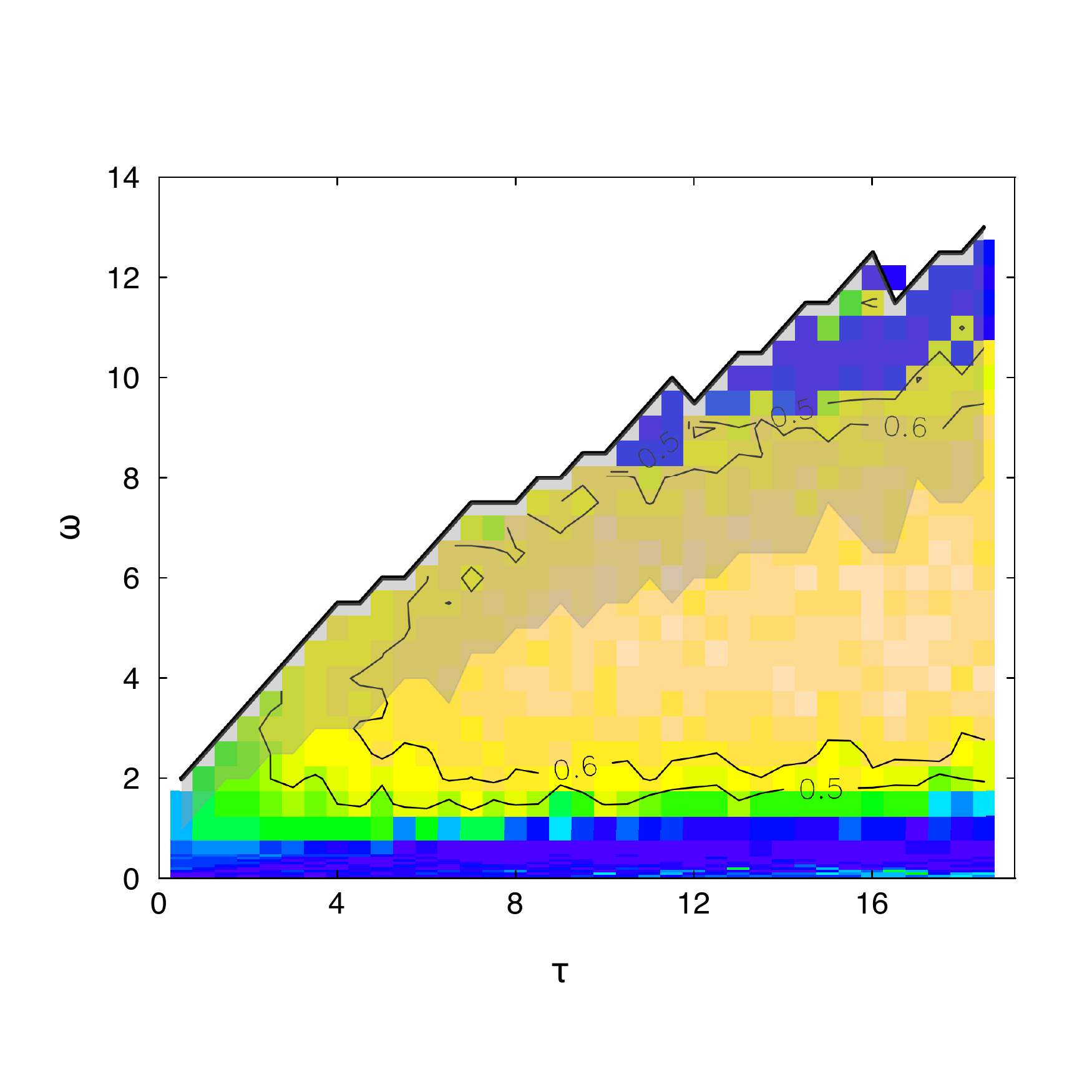}
	\caption{Equivalent of the bifurcation diagram for the stochastic model in the $(\tau,\omega_{SI}$) parameter space for $N=200$, $\gamma=1$,$\alpha_{SS}=0.04$. Identification of the oscillatory regime relies on the value of the frequency of peak power. Two potential boundaries are provided in the form of iso-lines at values $0.5$ and $0.6$. Peak frequency was $\approx 0.75$ (orange colour). Near zero frequencies are shown in dark blue. These boundaries are qualitatively consistent with those observed in the theoretical model. The thick black line shows one boundary for the disease-free regime determined as the value of $\omega_{SI}$ above which all realisations die out. The bottom boundary of the shaded area represents an alternative boundary determined as the value of $\omega_{SI}$ under which no realisations die out.}\label{fig:stochbifdiag}
\end{figure}

\textit{Boundary of the oscillating regime.} For those realizations that did not die out (over the entire duration of the run, i.e., $T_{max}=1320$), evidence for oscillatory behaviour was assessed through a
rigorous statistical analysis of the power spectrum of the time series of the number of infected nodes. The spectral estimation procedure was based on weighted periodogram estimates. The data was split into non overlapping segments, each containing $2^{13}$ points and the periodograms were calculated from application of discrete Fourier transforms (DFT). As application of DFT requires stationarity of
the time series, transients in the time series were removed using the following procedure. After confirming pseudo-stationarity of the last $2^{13}$ points (one segment) of the time series, their mean and
standard deviation was computed. A running average (low-pass filter) of the entire time series was then applied and the first time point after which all subsequent time points of the filtered time series
stayed within $5\%$ of the standard deviation of the mean of the last segment was selected as the starting point of the transient-free time series. Only those time series with at least $2^{15}$ time points
were kept for analysis, allowing for a theoretical minimum of $4$ segments to be included in the spectral estimation. In practice, provided transient-free data could be found, all parameter configuration yielded no less than $13$ segments with most configurations yielding $1300$ segments. Presence of a non-trivial oscillatory component in the signal was assessed by the presence of statistically significant power at a non-zero frequency peak in the power spectrum, with confidence intervals for the spectral estimates obtained as per the framework of~\cite{Halliday:1995ur}. The non-zero frequency constraint is required because, as mentioned in~\cite{Rogers2012}, time series in which the spectrum is dominated by the zero mode can be difficult to distinguish from pure white noise. To identify the boundaries of the oscillatory regime, power and frequency of the main peak of the power spectrum were recorded for all configuration pairs studied above. Figure~\ref{fig:oscboundary} shows a representative example (at $\tau=12$ but see Figure~\ref{fig:stochbifdiag} as well for an overall picture of the peak frequency) when $\omega$ is varied over the whole range $[0.05,13]$ and all other parameters are kept constant. The Figure reveals a pattern in which, with increasing $\omega$, there is initially a rapid increase in the power at the frequency of peak power, however, this frequency is small enough that the signal may not be considered oscillatory. This phase is followed by a clear increase in the frequency at peak power followed by a plateau with near constant power until $\omega$ reaches the boundary of the disease-free regime. As illustrated by Figure~\ref{fig:stochbifdiag}, with increasing $\tau$, the peak frequency can decrease to near zero levels such that the system briefly returns to the endemic regime. This is consistent with the result from the theoretical bifurcation diagram which shows a narrow strip of endemic regime between oscillatory regime and disease-free regime (e.g., see $\tau=4$ and $\omega\approx 780$ in Figure~\ref{SMF_Bif}). Thresholding of the frequency at peak power makes it possible to define a boundary for the oscillatory regime. This is a soft boundary in that there is no {\it a priori} basis for deciding a set level.

\textit{Bifurcation diagram in the $(\tau , \omega)$ plane.} Figure~\ref{fig:stochbifdiag} shows the bifurcation ``curves'' and the domains of the three different behaviours. Qualitatively, the plot confirms the prediction of the theoretical model with a bounded oscillatory domain abbuting the disease-free domain (here represented by a zero-disease boundary defined by all realizations dying out before $T_{max}$) at low values of $\tau$ ($\tau<10$) and showing a narrow strip of remaining endemic regime at higher $\tau$ values ($\tau>10$). When considering the stricter boundary, one observes an overlap between the disease-free regime and the oscillatory regime. Practically, this is characterised by time series of the number of infected nodes showing a large seesaw pattern of a long continuous decrease followed by a very rapid and large increase.

\subsection{Comparison of pairwise and simulation models}

\subsubsection{Comparison of the time evolution}

\begin{figure}[h]
  \centering
	\includegraphics[scale=0.35]{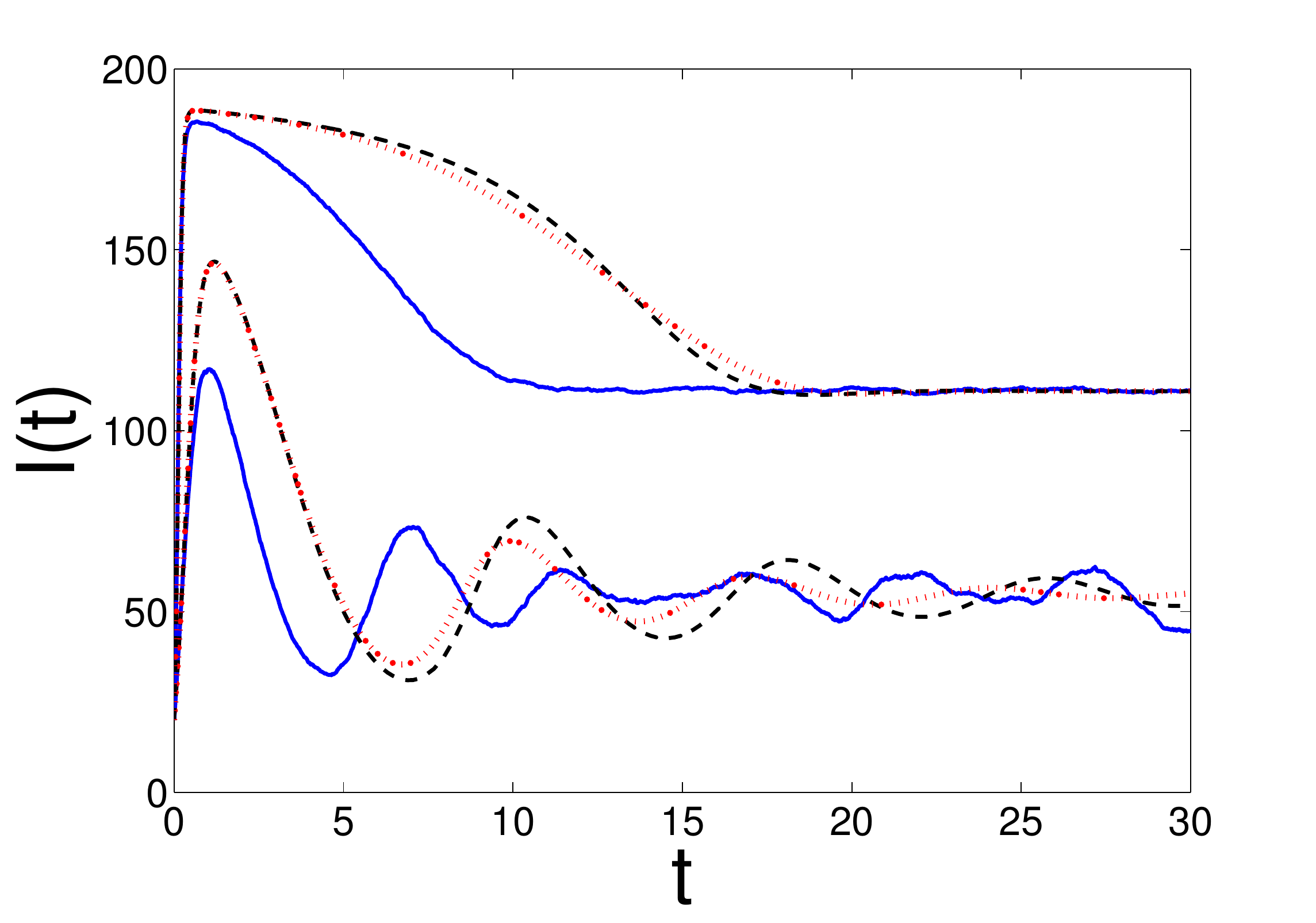}
	\caption{Comparison of the time evolution of prevalence in a network with $N=200$ nodes obtained from the average of 500 simulations (continuous blue curves), from the simple pairwise model \eqref{SMF1}-\eqref{SMF4} (dashed black curves) and from the compact pairwise model \eqref{CPw1}-\eqref{CPw5} (dotted red curves). The parameter values for the endemic case (upper curves) are $\alpha_{SS}=0.01$, $\omega_{SI}=0.7$, $\gamma=1$, $\tau=2$ and those yielding the oscillating solutions (lower curves) are $\alpha_{SS}=0.01$, $\omega_{SI}=1.5$, $\gamma=1$, $\tau=0.8$. For the oscillating case only simulations which did not die out until $t=20$ were taken into account.} \label{fig:comparison}
\end{figure}

In Fig.~\ref{fig:comparison}, we show direct comparisons between the outcome from stochastic simulation and results from the pairwise model. We illustrate the two typical behaviours in the presence of infection. As shown by the figure, there is reasonable agreement. Comparing oscillations is particularly difficult since small variations and de-phasing due to stochastic effects will lead to an average that smooths out the oscillatory behaviour. Hence, here we plot single realisations. Typically, it is relatively easy to find many parameter combinations where the agreement is satisfactory to good, especially for the endemic case and the endemic equilibrium. However, in what follows we aim to make a more principled comparison between the two models in the full parameter space.

\subsubsection{Comparison between bifurcation diagrams}

From a visual inspection of the two bifurcation diagrams in Figs.~\ref{SMF_Bif} and \ref{fig:stochbifdiag}, it is clear that the overall qualitative features of the model are captured well by both full stochastic and pairwise models. However, there are significant quantitative discrepancies. Most remarkably, the pairwise model significantly overestimates the slope of the transcritical bifurcation curves, effectively prescribing much higher cutting rates than needed in order to curtail an epidemic in the stochastic model. Starting from this crucial observation we comment on some fundamental and perhaps unresolvable differences which are mainly due to the natural limitations of the models. In particular, the transcritical bifurcation curve from the pairwise model requires that the network be fully connected. Hence, this naturally leads to a very high critical $\omega_{SI}$ rate, which is required in order to stop the epidemic in a fully connected network. There is a subtle point to be made here. Controlling the epidemic on a fully connected network via link cutting will, in most cases, require a punishingly high rate. This can take the network's average degree to very low values. This does not create a problem for the pairwise model in which the infection will be able to recover from very small levels of prevalence and spare networks. However, in such cases, in the stochastic model, there is a high probability that at low prevalence and sparsely connected network, the epidemic dies out. This suggests that smaller values of cutting rates are sufficient to yield the disease-free regime.

\subsubsection{Compact pairwise model}

Another important issue that may lead to disagreement is the heterogeneity of the degree distribution that is created by the process itself, and this will be studied in more detail in Section \ref{section_osc1}. Epidemics on networks with heterogeneous degrees can be described by heterogeneous mean field models \cite{EamesKeeling2002}. However, the number of equations in this type of models is of order $O(N^2)$ since the degree in a dynamic network can, in principle, vary between 0 and $N-1$, where $N$ is the number of nodes in the network. As a compromise between keeping degree heterogeneity and having a tractable system of ODEs, so-called compact pairwise models have been introduced \cite{HouseUnifying}. The variables of this model are $[S_k]$ and $[I_k]$ representing the average number of susceptible and infected nodes of degree $k$, respectively, and,  the average number of pairs $[SI]$, $[SS]$ and $[II]$. The significant reduction in the number of equations is achieved or made possible by not having to account for pairs like $[S_kI_l]$. Their precise bookkeeping is replaced by the following approximations
\begin{equation}
[A_kB] = [AB] \frac{k[A_k]}{\sum j[A_j]} , \label{pair_approx}
\end{equation}
where $[A_kB]=\sum_j [A_kB_j]$. In fact, the pairs $[A_kB_j]$ are not needed in the compact pairwise model, and only the pairs of the form $[A_kB]$ are used. We extended this model with terms responsible for link creation and deletion as follows. The deletion of links connecting an infected to a susceptible node with degree $k$ at rate $\omega_{SI}$ contributes positively to $[S_{k-1}]$ and negatively to $[S_{k}]$. The creation of links connecting a susceptible to another susceptible node with degree $k$ contributes negatively to $[S_{k}]$ with rate $\alpha_{SS}([S_{k}]([S]-1)-[S_{k}S])$ because the total number of such possible links is $[S_{k}]([S]-1)$ and the number of existing links is $[S_{k}S]$. The same process contributes positively to $[S_{k+1}]$. Using similar arguments we arrive to the following system,
\begin{align}
\dot{[S_k]}&=-\tau [S_kI] + \gamma [I_k]+ \omega_{SI}([S_{k+1}I]-[S_kI]) \label{CPw1}\\
 &+\alpha_{SS} ([S_{k-1}]([S]-1)-[S_{k-1}S])-\alpha_{SS}([S_{k}]([S]-1)-[S_{k}S])  , \nonumber  \\
\dot{[I_k]}&= \tau [S_kI] - \gamma [I_k]+\omega_{SI}([I_{k+1}S]-[I_kS]),  \label{CPw2} \\
\dot{[SI]}& = \gamma ( [II] - [SI])+ \tau ([SSI]-[ISI]-[SI])-\omega_{SI}[SI],  \label{CPw3} \\
\dot{[II]}& =-2\gamma [II] +2\tau ([ISI]+[SI]),  \label{CPw4} \\
\dot{[SS]}& =2\gamma [SI] - 2\tau [SSI]+\alpha_{SS} ([S]([S]-1)-[SS]),  \label{CPw5}
\end{align}
where $[S]=\sum_k [S_k]$, the approximation in Eq.~\eqref{pair_approx} is used to to compute the pairs $[S_kI]$ and $[I_kS]$ for $k=0,1, \ldots , N-1$, and according to \cite{HouseUnifying}, the triples are closed by
$$
[ASI]= \frac{[AS][SI]}{([SS]+[SI])^2}\sum_k  k(k-1)[S_k] .
$$
The solution of this system is also compared to simulation, see Fig.~\ref{fig:comparison}. One can observe that introducing degree heterogeneity in the model does not improve the agreement significantly. This may be explained by the fact that the main argument for closing triples in the pairwise and compact pairwise models is identical. Namely, both models assume that the states of a susceptible node's neighbours are effectively chosen independently and at random from the pool of available nodes. Nevertheless, the compact pairwise model may have relevance in studying the process because it enables us to investigate how the degree distribution varies in time, and this can become valuable especially in the oscillatory regime. The detailed study of the compact pairwise model is beyond the scope of the present paper, but it will certainly motivate further research.

\section{Explanation for the oscillatory behaviour}  \label{sec:ExplOsc}

\subsection{Simulation-based tracking}  \label{section_osc1}
Oscillations within the context of adaptive networks have proved to be difficult to map out, especially directly from simulations.
For many oscillatory systems the basic ideas of processes leading to oscillations are relatively simple. Usually, a combination of a positive and negative feedback with a suitable time delay leads to robust oscillations. Taking inspiration from this idea, it is possible to give a heuristic explanation for the appearance of oscillations in this adaptive network. The basic oscillating quantities in our explanation are the prevalence $[I]$ and the average degree $\langle k \rangle$, for which exact differential equations can be written down
\begin{align*}
\dot{[I]}&= \tau [SI] - \gamma [I],  \\
\langle \dot{k} \rangle &= \alpha_{SS} ([S]([S]-1)-[SI])-2\omega_{SI}[SI] .
\end{align*}
Thus as it can be immediately seen, the direction of change of these quantities is determined by the velocity of the following four processes.
\begin{itemize}
\item \textbf{A}: Infection with rate $\tau [SI]$,
\item \textbf{B:} Recovery of $I$ nodes with rate $\gamma [I]$,
\item \textbf{C:} Creation of $SS$ links with rate $\alpha_{SS}([S]([S]-1)-[SS])$ and
\item \textbf{D:} $SI$ link deletion with rate $\omega_{SI}[SI]$.
\end{itemize}
The important stages during the cycle of an oscillation are determined by the strength of process \textbf{A} relative to \textbf{B} and that of \textbf{C} compared to \textbf{D}. To explain and visualise this, let us consider an epidemic that is well established and is capable of sustaining oscillations.
Let us start from the situation where there are few infected nodes and the epidemic is about to take off, i.e., when process \textbf{A} dominates \textbf{B} and \textbf{C} is stronger than \textbf{D}. In this case, the expansion of the epidemic usually is also followed by
the network becoming more connected given that the total rate of link cutting is smaller than the total rate of link creation, i.e., $2\omega_{SI}[SI] < \alpha_{SS}([S]([S]-1)-[SS)$.
This is expected as the epidemic is just recovering from an excursion where the connectivity of the network was low and the number of susceptible nodes was large.
However, as the epidemic grows the balance of processes \textbf{C} and \textbf{D} changes. Namely, when the epidemic is still strong, the number of susceptible nodes decreases while the link cutting acts on an increasing number of ($SI$) edges. This reverts the balance of edge cutting and deletion meaning that now  $2\omega_{SI}[SI] > \alpha_{SS}([S]([S]-1)-[SS)$.  In Fig.~8a, a snapshot at the peak connectivity is shown. At this stage, the epidemic can and will continue to grow, i.e., \textbf{A} still dominates \textbf{B}. However, close to the highest possible prevalence level, shown in Fig.~8b, the recovery acting on the majority of the nodes and the cutting of occasional $(SI)$ links will outcompete the spread of infection, that is, \textbf{B} will be stronger than \textbf{A}. The continued loss of links leads to observing the `segregation' of the network, as shown in Fig.~8c, which amounts to the formation of susceptible and infected node clusters/clumps with occasional between cluster links. At this point, due to the increasing number of susceptible nodes, link creation will start to dominate and the susceptible clusters/clumps will become more densely connected. The survival of the epidemic now simply relies on the few inter-clump links which will allow infection to become re-established in the densely connected susceptible parts of the network. This is often followed by a sudden epidemic expansion, i.e., the dominance of \textbf{A} over \textbf{B}, see Fig.~8d. This last phase finishes the cycle and the system is back to the growing epidemic stage again as explained at the start of our argument. Summarising, the four stages of the oscillation can be characterised as follows:

1. \textbf{A} $>$ \textbf{B}, \textbf{C} $>$ \textbf{D}, $[I]$ increasing, $\langle k \rangle$ increasing,

2. \textbf{A} $>$ \textbf{B}, \textbf{C} $<$ \textbf{D}, $[I]$ increasing, $\langle k \rangle$ decreasing,

3. \textbf{A} $<$ \textbf{B}, \textbf{C} $<$ \textbf{D}, $[I]$ decreasing, $\langle k \rangle$ decreasing,

4. \textbf{A} $<$ \textbf{B}, \textbf{C} $>$ \textbf{D}, $[I]$ decreasing, $\langle k \rangle$ increasing.

\subsection{Exact master equation}

In order to understand the behaviour of the system we consider the exact master equations. First, let us describe the state space. Let $N$ denote the number of nodes and $E=N(N-1)/2$ the maximum number of edges. Fixing a network, the total number of states is $2^N$ since each node can be either an $S$ or an $I$ node. Each edge can be switched on or switched off, hence the total number of networks on $N$ nodes is $2^E$. For each network we have $2^N$ states, hence there are $2^{N+E}$ states in the state space. Even for small values of $N$ the state space becomes extremely large. For the sake of clarity however, the simplest cases of $N=2$ and $N=3$, i.e., for dynamic graphs with two or three nodes, are illustrated in Appendix B. For the case $N=2$, the transition matrix $M$ is also determined. This matrix has $2^{N+E}$ rows and columns and contains the transition probabilities between all possible states.

The master equation, describing the full system behaviour, takes the form $\dot x(t) = M x(t)$. Considering the special cases of small networks one can observe that this transition matrix has a special structure, however, we do not study this structure in detail. Instead, we turn to the question of the presence of oscillations in this exact model.
In order to answer this question the spectrum of the matrix $M$ has to be studied. It is known that for any transition matrix the spectrum is in the left half of the complex plain and zero is an eigenvalue. The eigenvector corresponding to the zero eigenvalue has non-zero coordinates only at those system states where all nodes are susceptible, and these states form an absorbing set. In the case of $N=2$ this eigenvector is $(\omega_{SS}, 0,0,0, \alpha_{SS},0,0,0)^T$. In the absorbing state all nodes are susceptible and the network changes dynamically according to the values of $\omega_{SS}$ and $\alpha_{SS}$. The solution $x(t)$ of the master equation is a linear combination of functions of the form $\mbox{e}^{\lambda t}u$, where $\lambda$ is an eigenvalue and $u$ is an eigenvector of $M$. If $\lambda \neq 0$, then it has a negative real part, hence this function tends to zero as $t\to \infty$, that is, for large $t$ only the absorbing state will be observed. However, typically there are eigenvalues with real part very close to zero yielding long time behaviour different from the absorbing state. For static networks this eigenvalue is real, hence the long time behaviour corresponds to a so-called quasi steady state, i.e., the expected prevalence is maintained at a non-zero constant value for a relatively long time. For dynamic networks this eigenvalue may be complex leading to damped oscillations, where slow damping leads to sustained, long-time oscillations.

Let us investigate briefly what sort of eigenvalues may be responsible for these oscillations. If the eigenvalue is complex, i.e., $\lambda = -\mu + \mbox{i} \nu$, then the dominant term in the expression of the prevalence is of the form $\mbox{e}^{-\mu t} \sin(\nu t)$, together with a similar term containing cosine instead of sine. A simple calculation shows that the ratio of two consecutive maxima is $\exp(-\mu/\nu)$. This ratio is always less than one showing that the oscillation is always damped, however, when $\mu/\nu$ is small then this damping is getting smaller. The key question is how the magnitude or smallness of this ratio depends on the values of the parameters. Even though this question is far beyond the scope of this paper, we mention a general, classical result related to this problem. In the early forties of the last century Kolmogorov posed the following problem. If $M$ is an $n\times n$ transition matrix of a continuous time Markov chain, find the location of its spectrum in the complex plane. Dmitriev and Dynkin proved that the spectrum is in the cone given by two half lines starting at the origin and having an angle $\frac{\pi}{2}- \frac{\pi}{n}$ with the negative part of the real axis \cite{DmitrievDynkin}. This implies that $\frac{-\mu}{\nu } \geq \tan (\frac{\pi}{n})$. It is easy to show that there exists a matrix $M$ for which equality can be achieved, i.e., there is a matrix that is optimal from the view point of oscillations. This matrix has ones below the diagonal and in the right top corner. This corresponds to a Markov chain with transitions aligned around a perfect cycle, thus giving rise to an idealised, perfect oscillator. The spectrum of our matrix $M$ in the master equation is in a much narrower cone, meaning that it yields much more markedly damped oscillations.

Future work will seek to estimate the angle of the smallest cone containing the spectrum of $M$, as given by the dynamic network model. Intuitively, it seems that the bigger the angle is, the longer the cycles of the Markov chain are. We note that in the case of $N=2$, for appropriate parameter combination, the longest cycle in the Markov chain's full transition diagram, see Fig.~\ref{fig_transN2}, has length four but with many potential perturbations or opportunities to deviate from the cycle ( the cycle with length four is highlighted by red arrows and the deviating transitions are shown by blue arrows in the figure). Summarising, we conjecture that the cycles in the Markov Chain's transition diagram are responsible for the oscillations observed in the probabilities of the states and in the expected values of some chosen quantities. Moreover, we claim that in the case of epidemic propagation in dynamic networks the Markov Chain's full transition map must contain cycles, with longer cycles and fewer potential off-cycles movements leading to stronger and more sustained oscillations.

\section{Network bifurcation}   \label{sec:NetBif}

\begin{figure}[h!]
	\begin{centering}
    	\includegraphics[scale=0.25]{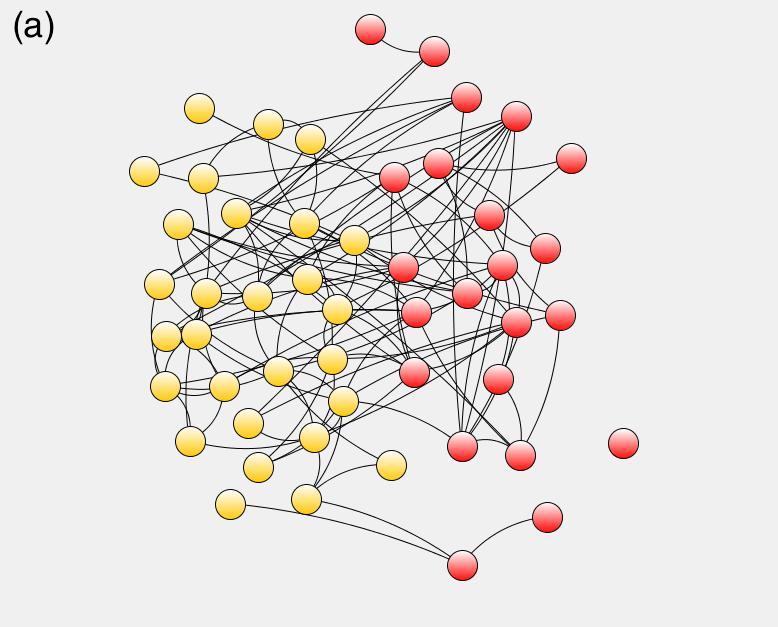}\label{fig:felmeno}
    	\includegraphics[scale=0.25]{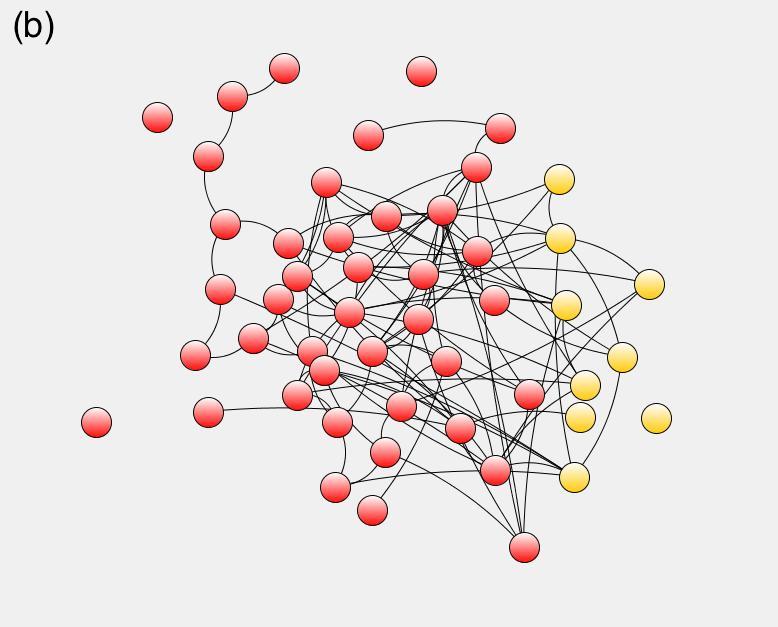}\label{fig:max}
	\end{centering}
	\newline
	
	\begin{centering}
	\includegraphics[scale=0.25]{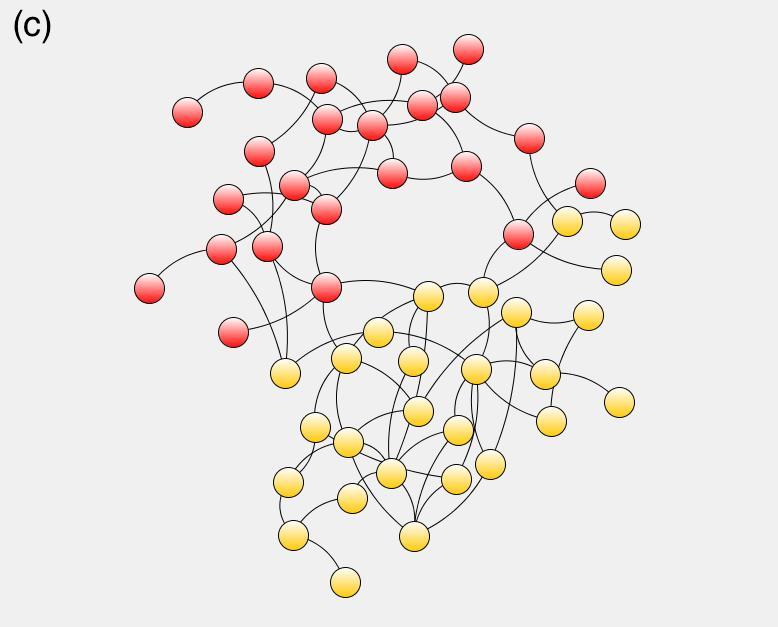}\label{fig:lenagy}
    	\includegraphics[scale=0.25]{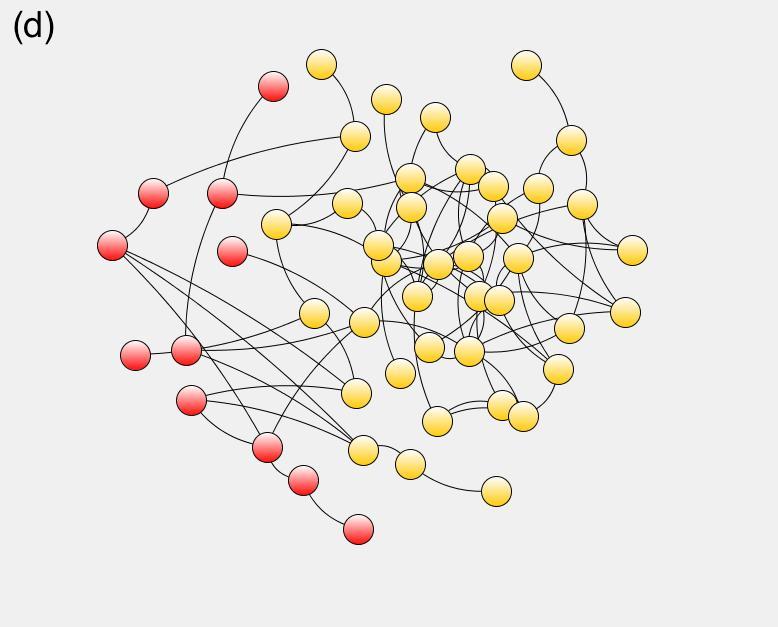}\label{fig:lekicsi}
		\end{centering}
\caption{Temporal snapshots of the main phases of an oscillation cycle are: (a) the growing phase of the epidemic with $\langle k \rangle$ close to its maximum, (b) close to the maximum prevalence and a decreasing average degree, (c) decreasing prevalence with $\langle k \rangle$ close to its minimum and, finally, (d) minimal prevalence but with growing average degree. Parameter values are $N=50$, $\tau=\gamma=1$, $\omega_{SI}=1.3$, $\alpha_{SS}=0.04$ with all the other activation and deletion rates being equal to zero.}
\label{fig:tempnetwshots}
\end{figure}

While some studies focusing on adaptive or dynamic networks consider and analyse changes in network structure,
a large proportion of papers either only or mainly focus on system level quantities such as infectious prevalence or some other population level
indicator with the main aim of characterising this via a bifurcation theory type analysis.
It has been observed that the networks themselves can also undergo significant changes in time or depending on parameters.
For example, the segregation of networks into different components has been observed by Gross et al. in \cite{ThiloPhyRevLett}, and see Fig.~\ref{fig:tempnetwshots}.
Such analysis can reveal important network features which can invalidate the use of mean-field or pairwise models and, more importantly, may reveal the impact of the
interplay between dynamics on and of the network. The emergence of network structure from such dynamic network models could be also interpreted
as a more natural or organic emergence of structure, as opposed to artificial or synthetic network models.
In what follows we aim to couple the analysis of system and network level changes, in order to
concurrently reveal the spectrum of behaviours at all levels.

In the case of scenario \textbf{A}, that is, when $\alpha_{SI}=\alpha_{II}=0$, $\alpha_{SS} \neq 0$, and $\omega_{II}=\omega_{SS}=0$, $\omega_{SI} \neq 0$, the most interesting types of networks or network properties that emerge from our model are: (a) networks with bimodal-like degree distribution and (b) networks where the largest connected component does not contain all nodes. These are all observed as transients and both features of the network manifest themselves at crucial or turning points in the epidemic dynamics. Namely, bimodal-like degree distributions are a feature of the network at low infection prevalence level when the susceptible and infected nodes segregate with few inter-component links. Here, link creation will increase the density in the susceptible cluster, while link cutting thins connections between infected nodes. Similarly, and also around this point, the clumping of similar states together, see Figs.~\ref{fig:lekicsi} and \ref{fig:lenagy}, leads to the network being composed of multiple disjoint components. Networks with these properties are transients. In the growing stages of the epidemic, as well as at the peak of it, networks are generally close to random or Erd{\H o}s-R\'enyi type networks.

\begin{figure}[!h]
\centering
\includegraphics[scale=0.5]{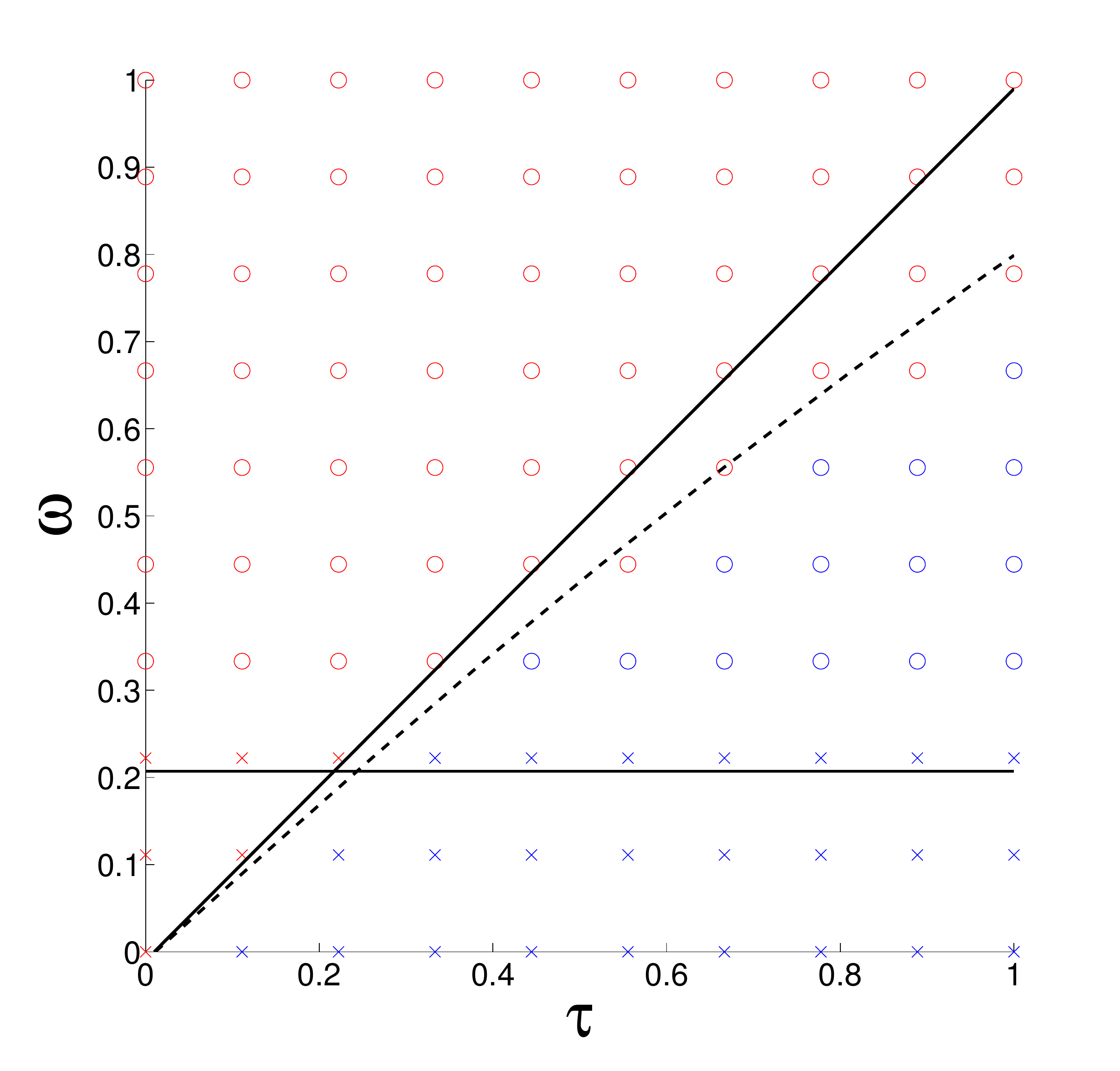}
\caption{Four different behaviours in the $(\tau, \omega)$ parameter space and the theoretical bifurcation curves for $\alpha = 0.01$, $N=200$, $\gamma=1$. The horizontal line represents the boundary of the parameter domain where the graph is connected. In the simulation, networks which on average had at least 3 disjointed components were considered disconnected. The other two curves are the transcritical bifurcation curves obtained from the mean-field approximation (continuous diagonal line) and from the pairwise approximation \eqref{tauTC} (dashed curve). The markers are as follows: \textcolor[rgb]{0.00,0.00,1.00}{$\times$} - connected, epidemic, \textcolor[rgb]{1.00,0.00,0.00}{$\times$} - connected, no epidemic, \textcolor[rgb]{0.00,0.00,1.00}{$\circ$} - disconnected, epidemic, and \textcolor[rgb]{1.00,0.00,0.00}{$\circ$} - disconnected, no epidemic.}
\label{fig:netwbif}
\end{figure}

In the case of scenario \textbf{B}, that is, when $\alpha_{SI}=\alpha_{II}= \alpha_{SS} :=\alpha$, and $\omega_{II}=\omega_{SS}=\omega_{SI} := \omega$ a more complete characterisation of network bifurcations can be achieved. As suggested by the pairwise model, two behaviours may occur according to the long time prevalence level, namely disease-free or endemic steady state. As regards the network structure we studied the connectivity  (through determining the type and number of connected components) of the network and the degree distribution. Our goal here was to map out system and network behaviour over the $(\tau, \omega )$ parameter space. For the points of a lattice in the $(\tau, \omega)$ parameter plane several stochastic simulations were run. The average epidemic level end network connectivity was then determined in the steady state (after a sufficiently long time). The different system and network level outcomes yielded four different combinations of behaviours as shown in Fig~\ref{fig:netwbif}. The most interesting observation is that epidemics can be curtailed either due to the network being disconnected or due to the epidemic being sub-threshold when the network could theoretically support an epidemic.


The bifurcation curves can also be determined by theoretical considerations. The results of these are shown by lines in the figure. The simplest way of finding the transcritical bifurcation curve, which separates the endemic and disease-free regions, is to start from the steady state equation $\tau [SI] = \gamma [I]$ and use the pair closure $[SI]\approx \frac{n}{N-1}[S][I]$, where $n$ denotes the average degree. It was shown in \cite{PRAcikk} that the steady state value of the average degree can be given as $n=\frac{\alpha(N-1)}{\alpha+\omega}$. Substituting this expression in the steady state equation, then dividing by $[I]$ and substituting $[S]=N$, which holds at the boundary of the endemic region, we arrive at
$$\tau_{pc} = \gamma \frac{\alpha+\omega}{N \alpha}.$$  This theoretical bifurcation value, based on the pair closure, is shown in Fig.~\ref{fig:netwbif} as a continuous diagonal line.
The transcritical bifurcation curve can also be determined from the pairwise model that is closed at the level of triples. As shown in Section \ref{sec:PW}, the transcritical bifurcation in that model occurs at the value of $\tau$ given by \eqref{tauTC}. Plotting this critical value as a function of $\omega$ we obtain the bifurcation curve also shown in Fig.~\ref{fig:netwbif} (dashed line).
It is worth noting that closing at the level of triples leads to a better agreement with simulation results, as can be seen in Fig.~\ref{fig:netwbif}.

The bifurcation curve separating the connected and disconnected region can be derived analytically as follows. In \cite{PRAcikk} it was shown that the network becomes an Erd\H os-R\'enyi type random graph at the steady state and the probability of edge activation is $p=\frac{\alpha}{\alpha + \omega}$. Moreover, we know that the threshold for an Erd\H os-R\'enyi graph being disconnected is $p=\frac{\ln N}{N}$ \cite{ChungLu,DurettRanGraph}, where $N$ denotes the size of the graph. Taking into account these two formulas we get the following equation,
$$p=\frac{\alpha}{\alpha + \omega}=\frac{\ln N}{N}.$$
Thus the critical threshold for connectivity is,
$$\omega^*=\alpha(\frac{N}{\ln N} - 1).$$
The horizontal line in Fig.~\ref{fig:netwbif} is drawn at this value of $\omega$.
One can see that this is in good agreement with the connectivity results obtained from simulations.
Moreover, the distribution of networks during and at the end of simulations is well described by the binomial distribution (figures not shown) and for fixed values of $\alpha$, the number of components increases sharply with larger values of $\omega$.

\section{Discussion}

We proposed a dynamic network model using $SIS$ epidemic propagation with both random and preferential, i.e., link-type based, link activation and deletion. We analysed the pairwise model and
carried out network-based stochastic simulations. The oscillations developing in the system were studied using Fourier analysis, as well as through analysing the exact stochastic model in terms of
the master equations. Further, we developed a compact pairwise approximation for the current dynamic network model and carried out a detailed study of the agreement between the ODE approximation
and simulations. For the case of random activation deletion we provided a bifurcation analysis of the network structure itself.

Both simple pairwise and compact pairwise models fail to accurately describe the oscillatory regime. One of the main reasons of this failure is that the closures in these systems are based on the
assumption that the distribution of infected nodes around susceptible nodes is binomial, which does not hold in simulation, especially close to the die-out regime. Fig.~8c, where infected and
susceptible clumps are developed, clearly shows that the binomial assumptions is completely inaccurate. The poor performance of pairwise models for adaptive networks is known in the literature, see
for example \cite{ThiloPhysD2014}. Further investigation of the time dependence of the network structure and the failure of the binomial assumption may reveal the reasons behind this phenomenon.

Further important questions that we wish to highlight are related to characterising and measuring changes in the structure or bifurcation of networks.
In this paper we focused mainly on scenario \textbf{B}, where we identified the connectedness of the network as an important indicator
of changes in network structure. This to some extent was dictated by the particular choice of model, where the epidemic and network dynamics are decoupled.
This choice is not unique and investigating the degree distribution or its average could be equally relevant or even more informative. For example, for scenario \textbf{B},
apart from considering the evolution of the average degree, one could for example map out the evolution of clustering and the evolution of the degree distribution which
undergoes interesting changes from binomial-like degree distributions to markedly bimodal degree distributions, when the networks is close to the segregation point.
Such scenarios could be further investigated using the compact pairwise in order to shed some light on the evolution of the degree distribution, 
and could be again used as a precursor to more detailed analysis based on simulations. Hence, there are many possibilities and the right choice may not be universal but model dependent. Nevertheless, we wish to emphasise the importance of clearly describing changes in network structure, which can then lead to a better understanding of how networks with certain properties emerge and
why mean-field models break down and what their limitations are.

Studying small toy networks and in particular the structure of the transition diagram of full Markov Chains, has allowed us to relate the oscillatory phenomena to cycles and their length in the Markov
Chains.  From our simple experiments, it is worth noting that the creation of $II$ links, i.e., non-zero value of $\alpha_{II}$, may lead to more pronounced oscillations, as this can lead to the appearance of
longer cycles in the transition diagram. Finally, determining the transition diagram of the lumped systems for small graphs and investigating the maximal length of cycles in the transition portrait of the Markov Chain and their
relation to the oscillations produced by the system could shed further light behind how oscillations emerge and can be sustained.

\section*{Acknowledgements} P\'eter L. Simon acknowledges support from OTKA (grant no. 81403).

\section*{Appendix A: Bifurcation analysis of the mean-field model} \label{AppendixA}

In this appendix the case of \textbf{Scenario A}, that is $\alpha_{SI}=\alpha_{II}=0$, $\alpha_{SS} \neq 0$, and $\omega_{II}=\omega_{SS}=0$, $\omega_{SI} \neq 0$ is considered.

\subsection*{A1: Endemic steady state} \label{SSendemic}

The steady states are given by system \eqref{SMF1}-\eqref{SMF4} by putting zeros in the left hand sides (and omitting those parameters that are assumed to be zero), i.e. by system
\begin{align}
0 &=\tau [SI] - \gamma [I],  \label{MF1} \\
0 &=\gamma ( [II] - [SI])+ \tau ([SSI]-[ISI]-[SI])-\omega_{SI}[SI],  \label{MF2} \\
0 &=-2\gamma [II] +2\tau ([ISI]+[SI]),  \label{MF3} \\
0 &=2\gamma [SI] - 2\tau [SSI]+\alpha_{SS} ([S]([S]-1)-[SI])  \label{MF4} \\
\end{align}
where the closures
\begin{equation}
[SSI]=\frac{(k_{S}-1)[SS][SI]}{k_S[S]} \,\,\, \text{and} \,\,\, [ISI]=\frac{(k_{S}-1)[SI][SI]}{k_S[S]} \label{clos}
\end{equation}
are used with
$$
k_{S}=\frac{[SS]+[SI]}{[S]} .
$$
For simplicity, the steady state values of the four variables are denoted here by $[I]$, $[SI]$, $[II]$ and $[SS]$. All of them will be expressed in terms of $[S]$ and then a quartic equation for $[S]$ will be derived. It is obvious that $[I]=N-[S]$, then from \eqref{MF1} we get $[SI]=\frac{\gamma}{\tau}[I]$. Multiplying \eqref{MF2} by 2 and adding it to \eqref{MF3} and to  \eqref{MF4} one can express $[SS]$ in terms of $[I]$ and $[S]$ as
\begin{equation}
[SS]=[S]([S]-1)-2\frac{\omega_{SI}\gamma}{\alpha_{SS}\tau} [I].
\label{SS}
\end{equation}
From \eqref{MF3} we can express $[II]$ in terms of $[I]$ and $[S]$ as
\begin{equation}
[II]=\frac{\tau [SI]}{\gamma} \left(\frac{[SI]}{[S]} - \frac{[SI]}{[SS]+[SI]} +1  \right).
\label{II}
\end{equation}
Multiplying \eqref{MF2} by 2 and adding it to \eqref{MF3} we get
$$
2\tau [SSI] -2\gamma [SI] -2\omega_{SI} [SI] = 0 .
$$
Using the closure and dividing by $2[SI]$ yields
$$
\frac{[SS]}{[S]} - \frac{[SS]}{[SS]+[SI]} = \frac{\gamma + \omega_{SI} }{\tau} .
$$
Now substituting $[SI]=\frac{\gamma}{\tau}(N-[S])$ and \eqref{SS} into this equation one obtains the following quartic equation for $x=[S]$
\begin{equation}
x^4+A_3 x^3 + A_2 x^2 + A_1 x + A_0 =0
\label{quartic}
\end{equation}
with
\begin{align*}
A_3 &= 4ab-3-2b-c\\
A_2 &= 2+2b+c+b^2+bc-6ab-4ab^2 -2abc + 4a^2b^2+Nb(1-4a)\\
A_1 &= Nb(-1+6a-b-c+6ab+2ac-8a^2b)\\
A_0 &= 2N^2ab^2(1-2a), \\
\end{align*}
where $a=\frac{\omega_{SI}}{\alpha_{SS}}$, $b=\frac{\gamma}{\tau}$ and $c=\frac{\omega_{SI}}{\tau}$ .

Once \eqref{quartic} is solved for $x$, then the steady state $([I], [SI], [II], [SS])$ is given by $[I]=N-x$, $[SI]=\frac{\gamma}{\tau}[I]$, \eqref{SS} and \eqref{II}. A solution $x$ yields a biologically meaningful steady state if all of its coordinates are non-negative. An extensive numerical study showed that for any combination of the parameters there can be at most one endemic steady state. The endemic steady state exists if the disease-free equilibrium is unstable. The condition for that is determined in the next subsection. We note that the above calculation yields only the endemic steady states because during the derivation of the quartic equation there was a division by $2[SI]$, which is zero at the disease-free steady state.

\subsection*{A2: Stability of the disease-free steady state}

The stability of the disease-free steady state  is determined by the Jacobian matrix determined at $([I], [SI], [II], [SS])=(0,0,0,N(N-1))$. In order to compute this matrix we determine the partial derivatives of the triples at this steady state, as they are given by the closures \eqref{clos}.
$$
\frac{\partial [SSI]}{\partial [I]} = 0  , \quad  \frac{\partial [SSI]}{\partial [SI]} = N-2 , \quad  \frac{\partial [SSI]}{\partial [II]} = 0 , \quad  \frac{\partial [SSI]}{\partial [SS]} =  0 ,
$$
$$
\frac{\partial [ISI]}{\partial [I]} = 0  , \quad  \frac{\partial [ISI]}{\partial [SI]} = 0 , \quad  \frac{\partial [ISI]}{\partial [II]} = 0 , \quad  \frac{\partial [ISI]}{\partial [SS]} =  0 ,
$$
Using these partial derivatives the Jacobian matrix at the disease-free steady state can be given as
$$
J =
 \begin{pmatrix}
  -\gamma & \tau & 0 & 0 \\
  0 & -\gamma + \tau (N-2) - \tau -\omega_{SI}) & \gamma & 0 \\
  0  & 2\tau  & -2\gamma & 0  \\
  \alpha_{SS}(1-2N) & 2\gamma-2\tau(N-2) & 0 & -\alpha_{SS}
 \end{pmatrix}.
$$
It can be easily seen that $-\alpha_{SS}$ and $-\gamma$ are eigenvalues of this matrix. Te remaining two eigenvalues are the eigenvalues of the $2\times 2$ matrix in the middle:
$$
\begin{pmatrix}

   -\gamma + \tau (N-2) - \tau -\omega_{SI}) & \gamma  \\
   2\tau  & -2\gamma
\end{pmatrix}.
$$
The disease-free steady state is stable if and only if all the eigenvalues have negative real part. This has to be checked only for the above $2\times 2$ matrix. For a $2\times 2$ matrix the eigenvalues have negative real part if and only if its determinant is positive and its trace is negative. The determinant is positive if $\gamma + \omega_{SI} -\tau (N-2)>0$. The trace is positive if $3\gamma + \omega_{SI} -\tau (N-3)>0$. The first condition implies the second one, hence we proved Proposition \ref{prop1}.

\subsection*{A3: Stability of the endemic steady state}

The stability of the endemic steady state can be determined only numerically. For a given set of the parameters the coordinates of the endemic steady state can be computed according to Appendix A1. The partial derivatives in the Jacobian $J$ can be calculated analytically, then substituting the numerically obtained coordinates of the endemic steady state we get the entries of the Jacobian numerically. This enables us to calculate the coefficients of the characteristic polynomial
\[
\lambda^4-b_3\lambda^3 +b_2\lambda^2-b_1\lambda +b_0,
\]
where $b_3 = \Tr J$, $b_0 = \det J$ and $b_1,b_2$ can be given as the sum of some subdeterminants of the Jacobian, the concrete form of which is not important at this moment. To find the parameter values where Hopf bifurcation occurs we use the method introduced in \cite{Szabo2012DEA}. In the case of $4\times 4$ matrices the necessary and sufficient condition for the existence of pure imaginary
eigenvalues is
\begin{equation} \label{hopf}
b_0b^2_3= b_1(b_2b_3-b_1)\text{ and }\sign b_1 = \sign b_3,
\end{equation}
Thus the Hopf-bifurcation set in the $(\tau,\omega_{SI})$ parameter plane can be obtained as follows. For a given value of $\tau$ we compute the value of $b_0b^2_3= b_1(b_2b_3-b_1)$ numerically as $\omega_{SI}$ is varied. It turns out that for a range of $\tau$ values this expression changes sign twice as $\omega_{SI}$ is varied. More precisely, for given values of the other parameters ($N, \gamma, \alpha_{SS}$) there exist $\tau$ values $\tau_1$ and $\tau_2$, such that for $\tau\in (\tau_1, \tau_2)$ we get $\omega_1$ and $\omega_2$ such that for $\omega_{SI}=\omega_i$ ($i=1,2$) we have $b_0b^2_3= b_1(b_2b_3-b_1)$, i.e. Hopf bifurcation occurs at $\omega_{SI}=\omega_i$ ($i=1,2$). If $\tau$ is not in the interval $(\tau_1, \tau_2)$, then there is no Hopf bifurcation, i.e. the relation $b_0b^2_3= b_1(b_2b_3-b_1)$ cannot hold. For $\tau\in (\tau_1, \tau_2)$ and $\omega_{SI}\in (\omega_1,\omega_2)$ there is a stable periodic orbit. If $\omega_{SI}$ is outside the interval $(\omega_1,\omega_2)$, then there is no periodic orbit and either the endemic or the disease-free steady state is stable. The final state of the system is shown in the bifurcation diagram in Fig.~\ref{SMF_Bif}.

\section*{Appendix B: State space and transition matrix for adaptive graphs with $N=2$ and $N=3$} \label{AppendixB}

\begin{figure}[!h]
\centering
\includegraphics[scale=0.8]{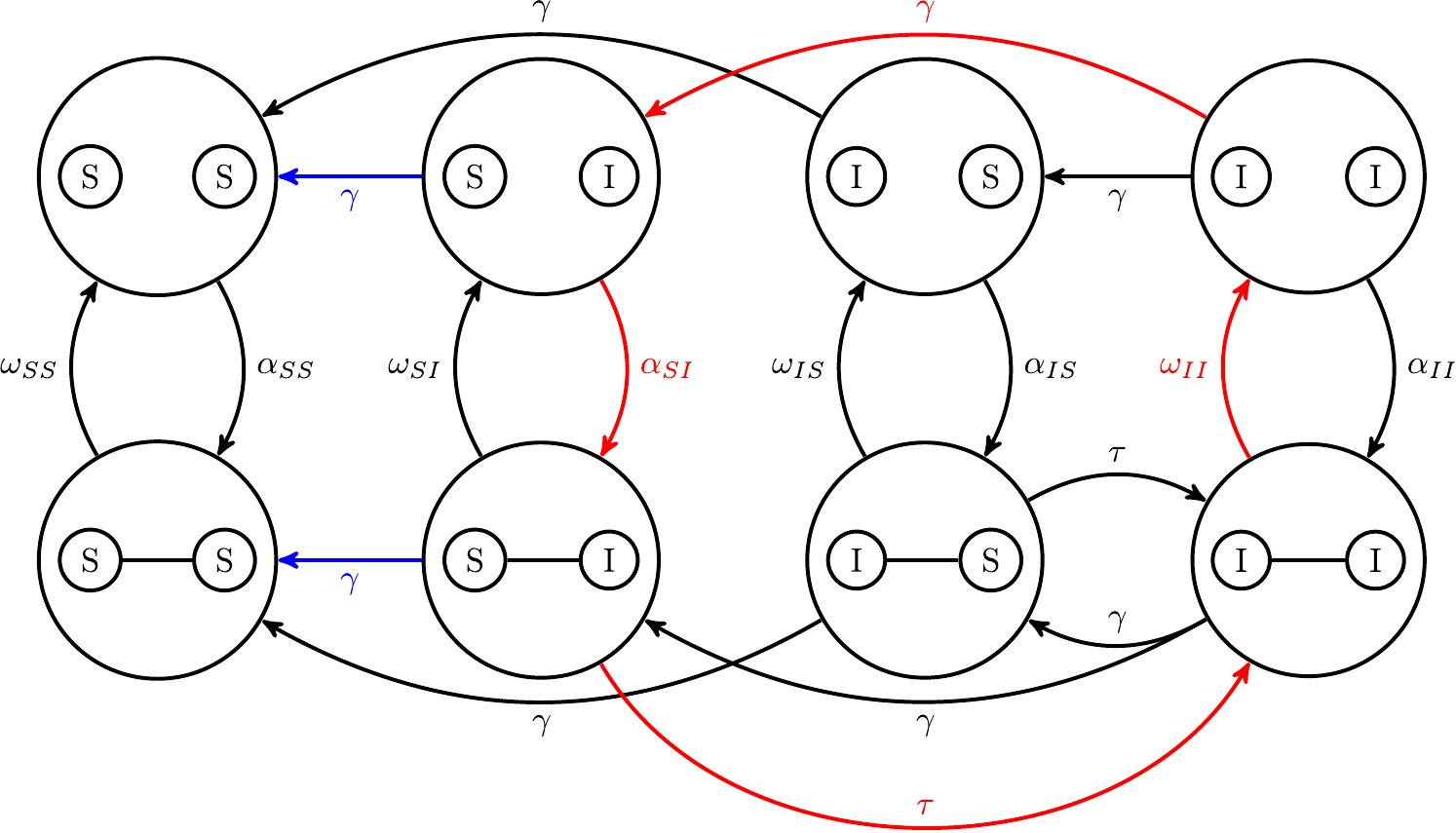}
\caption{The state space and the transitions for a dynamic network with two nodes. Red arrows represent the longest possible cycle of length four, with the blue arrows highlighting transitions
that would drive the process off cycle.} \label{fig_transN2}
\end{figure}

In the case $N=2$, i.e. for a graph with two nodes, the number of edges is at most $E=1$. That is there are two graphs on two nodes, one is a single edge the other one consists of two disjoint nodes. We denote the states with $\{ SS, SI, IS, II\}$ when the graph consists of two disjoint nodes. Here the state $SI$ means that node 1 is $S$ and node 2 is $I$. The states are denoted by $\{ \overline{SS}, \overline{SI}, \overline{IS}, \overline{II} \}$ when the graph is a single edge. Thus the full state space for $N=2$ contains the following 8 states:  $\{ SS, SI, IS, II, \overline{SS}, \overline{SI}, \overline{IS}, \overline{II}\}$. Consider now the transitions between these states. There are two types of transitions: epidemic transitions (infection and recovery) and network transitions (creating and deleting edges). Epidemic transitions may occur among states that belong to the same graph, that is within the subsets $\{ SS, SI, IS, II\}$ and $\{ \overline{SS}, \overline{SI}, \overline{IS}, \overline{II} \}$. Within the first subset only recovery may occur since these states belong to a graph that consists of two separate nodes. So the only possible transitions are $SI\to SS$, $IS\to SS$, $II\to IS$ and $II\to SI$, these may happen with rate $\gamma$. Within the subset $\{ \overline{SS}, \overline{SI}, \overline{IS}, \overline{II} \}$ infection may happen as well, hence the possible transitions are $\overline{SI}\to \overline{II}$ and $\overline{IS}\to \overline{II}$ with rate $\tau$ and $\overline{SI}\to \overline{SS}$, $\overline{IS}\to \overline{SS}$, $\overline{II}\to \overline{IS}$  and $\overline{II}\to \overline{SI}$ with rate $\gamma$. Network transitions occur between states in which the corresponding nodes are of the same type. For example, the transition $SS\to \overline{SS}$ occurs at rate $\alpha_{SS}$ since an $SS$ type edge is created during this transition. Similarly, the transition $\overline{SS}\to SS$ occurs at rate $\omega_{SS}$ since an $SS$ type edge is deleted during this transition. In general, the transition $AB\to \overline{AB}$ happens at rate $\alpha_{AB}$ and the transition $\overline{AB}\to AB$ happens at rate $\omega_{AB}$, where $A,B\in \{S,I\}$. All the transitions are shown in Figure \ref{fig_transN2}. If the states are ordered as $\{ SS, SI, IS, II, \overline{SS}, \overline{SI}, \overline{IS}, \overline{II}\}$, then the transition matrix of the corresponding Markov chain takes the form
$$
M= \left(
\begin{array}{cccc}
M_{11}-A & \Omega\\
A & M_{22}- \Omega
\end{array}
\right),
$$
where
$$
M_{11}= \left(
\begin{array}{cccc}
0 & \gamma &  \gamma  & 0\\
0 &-\gamma & 0 &  \gamma \\
0 & 0 &-\gamma & \gamma  \\
0 & 0 & 0 & -2\gamma
\end{array}
\right), \quad
M_{22}= \left(
\begin{array}{cccc}
0 & \gamma & \gamma & 0\\
0 & -\tau -\gamma& 0 & \gamma \\
0 & 0 & -\tau -\gamma& \gamma \\
0 & \tau & \tau & -2\gamma
\end{array}
\right),
$$
$$
A = \left(
\begin{array}{cccc}
\alpha_{SS} & 0 &  0  & 0 \\
0 & \alpha_{SI} & 0 &  0 \\
0 & 0 & \alpha_{SI}  & 0 \\
0 & 0 & 0 & \alpha_{II}
\end{array}
\right), \quad
\Omega= \left(
\begin{array}{cccc}
\omega_{SS} & 0 & 0 & 0\\
 0 & \omega_{SI} & 0 & 0\\
 0 & 0 & \omega_{SI} & 0\\
 0 & 0 & 0 & \omega_{II}
 \end{array}
\right).
$$
The matrices $M_{11}$ and $M_{22}$ contain the transition rates corresponding to the epidemic transitions. These rates belong to the transitions within the subsets $\{ SS, SI, IS, II\}$ and $\{ \overline{SS}, \overline{SI}, \overline{IS}, \overline{II} \}$, respectively. The matrices $A$ and $\Omega$ contain the network transition rates that correspond to transitions between these two subsets. The master equation can be written as $\dot x(t) = M x(t)$, where the coordinates of the eight dimensional vector $x(t)$ are the probabilities of the states at time $t$.

Let us briefly consider the case $N=3$. Then $E=3$, hence there are $2^3=8$ different possible graphs, one without any edges, three graphs with one edge, three line graphs with two edges and a triangle with three edges. Each graph can be in 8 possible states, namely $\{ SSS, SSI, SIS, ISS, SII, ISI, IIS, III\}$. Hence there are $2^3\cdot 2^3=64$ states altogether. Epidemic transitions may occur among states that belong to the same graph, for example, in the case of a triangle graph the transition $IIS\to III$ happens at rate $2\tau$, while its rate is $\tau$ for a line graph where node 2 is connected to node 1 and node 3. The transition rate may be zero, e.g. in the case when there is only one edge in the graph connecting nodes 1 and 2, or in the case when there are no edges at all in the graph. Network transitions occur between states in which the corresponding nodes are of the same type. For example, denoting by $SSS$ the state where all nodes are susceptible and the graph consist of three disjoint nodes, and by $\overline{SS}S$ the state where all nodes are susceptible and the graph contains one edge that connects node 1 and 2, the rate of transition $SSS\to \overline{SS}S$ is $\alpha_{SS}$, while the rate of transition $\overline{SS}S\to SSS$ is $\omega_{SS}$. The master equations take again the form $\dot x(t) = M x(t)$, where the coordinates of the 64-dimensional vector $x(t)$ are the probabilities of the states at time $t$.

We note that the size of the state space can be reduced by lumping some states together, similarly to the case of static graphs \cite{SimonTaylorKiss}. The lumping of the state space for dynamic network processes is beyond the scope of this paper, here we only mention a few simple cases where lumping can be carried out. In the case $N=2$ it is easy to see that the states $SI$ and $IS$ can be lumped together, which means that their sum can be introduced as a new variable, and their differential equations can be added up. Similarly, the sum of $\overline{SI}$ and $\overline{IS}$ can be introduced as a new variable. (By adding their differential equations one can immediately see that the old variables will not appear in the remaining system of equations.) Hence the eight-dimensional system $\dot x(t) = M x(t)$ can be reduced to a six-dimensional system by lumping. In the case $N=3$ we have even more chance for lumping if it is assumed that $\alpha_{SI}=0$, $\alpha_{II}=0$, $\omega_{SS}=0$, $\omega_{II}=0$. Without explaining the details we claim that in this case the 64-dimensional system of master equations can be reduced to a 20-dimensional one by lumping. In the case $N=4$ the dimension of the state space can be reduced from 1024 to 89.

\newpage

\newpage

\newpage
%

\newpage

\newpage

\newpage

\newpage

\newpage

\newpage

\newpage

%

\newpage

\end{document}